\documentclass{article}

\usepackage{Arxiv}

\usepackage{MyNotations}

\usepackage{amssymb}

\usepackage[backend=biber,bibstyle=custom-apa,labeldate=year,citestyle=apa,isbn=false,url=false,eprint=false,backref=true,sortcites=false,sorting=nyt]{biblatex} 
\newcommand{\mkbibbracketscol}[1]{\mkbibbrackets{\textcolor{accentcolor!80}{#1}}}
\DeclareCiteCommand{\parencite}[\mkbibbracketscol]
{\usebibmacro{cite:init}%
	\usebibmacro{prenote}%
	\toggletrue{apa:inpcite}}
{\usebibmacro{citeindex}%
	\printtext[bibhyperref]{\usebibmacro{cite}%
		\usebibmacro{cite:post}%
		\togglefalse{apa:inpcite}}}
{\multicitedelim}
{\usebibmacro{postnote}}

\usepackage{microtype}
\usepackage{doi}

\newcommand*{\KS}{\textcolor{black}}
\newcommand*{\AD}{\textcolor{black}}

\usepackage[capitalise,noabbrev,nameinlink]{cleveref}

\title{ Finite Element Neural Network Interpolation. \\ Part I: Interpretable and Adaptive Discretization for Solving PDEs}

\author{
	Kateřina Škardová  \thanks{Corresponding author} \\
	LMS, École Polytechnique,\\
	IPP/CNRS, Palaiseau, France.\\
  	\texttt{} \\
	INRIA, Palaiseau, France.\\
	\href{mailto:katerina.skardova@polytechnique.edu}{\texttt{katerina.skardova@polytechnique.edu}} \\
	\And
	Alexandre Daby-Seesaram \\
	LMS, École Polytechnique,\\
	IPP/CNRS, Palaiseau, France.\\
	\texttt{} \\
	INRIA, Palaiseau, France.\\
	\href{mailto:alexandre.daby-seesaram@polytechnique.edu}{\texttt{alexandre.daby-seesaram@polytechnique.edu}} \\
	\And
	Martin Genet \\
	LMS, École Polytechnique,\\
	IPP/CNRS, Palaiseau, France.\\
  	\texttt{} \\
	INRIA, Palaiseau, France.\\
	\href{mailto:martin.genet@polytechnique.edu}{\texttt{martin.genet@polytechnique.edu}} \\
}

\usepackage{tocloft}  

\begin{document}

		\maketitle

		\begin{abstract}
			We present the Finite Element Neural Network Interpolation (FENNI) framework, a sparse neural network architecture extending previous work on Embedded Finite Element Neural Networks (EFENN) introduced with the Hierarchical Deep-learning Neural Networks (HiDeNN). Due to their mesh-based structure, EFENN requires significantly fewer trainable parameters than fully connected neural networks, with individual weights and biases having a clear interpretation.
			
			Our FENNI framework, within the EFENN framework, brings improvements to the HiDeNN approach. First, we propose a reference element-based architecture where shape functions are defined on a reference element, enabling variability in interpolation functions and straightforward use of Gaussian quadrature rules for evaluating the loss function. Second, we propose a pragmatic multigrid training strategy based on the framework's interpretability. Third, HiDeNN’s combined rh-adaptivity is extended from 1D to 2D, with a new Jacobian-based criterion for adding nodes combining h- and r-adaptivity. From a deep learning perspective, adaptive mesh behavior through rh-adaptivity and the multigrid approach correspond to transfer learning, enabling FENNI to optimize the network's architecture dynamically during training.
			
			The framework's capabilities are demonstrated on 1D and 2D test cases, where its accuracy and computational cost are compared against an analytical solution and a classical FEM solver. On these cases, the multigrid training strategy drastically improves the training stage’s efficiency and robustness. Finally, we introduce a variational loss within the EFENN framework, showing that it performs as well as energy-based losses and outperforms residual-based losses.
			
			This framework is extended to surrogate modeling over the parametric space in Part II.
		\end{abstract}


	\newpage

\section{Introduction}
Solving partial differential equations (PDEs) is essential for many real-world applications. Given that analytical solutions are rarely obtainable for more complex PDEs, a range of computational methods for numerical solutions of PDEs have been developed \parencite{levequeFiniteDifferenceMethods2007,brennerMathematicalTheoryFinite2008,cardiffThirtyYearsFinite2021,liuEightyYearsFinite2022}. \AD{These methods typically replace differential operators and the solution space with discrete approximations, transforming the problem into a set of finite-dimensional algebraic equations}

\AD{However, for applications requiring several computations or real-time simulations, those numerical methods may still be computationally too expensive. Reduced-order methods (ROM) address this issue by reducing the problem's dimensionality.} Several reduced-order modelling techniques rely on solving full-order computations for specific sets of parameters \parencite{sirovich_turbulence_1987,chatterjee_introduction_2000,maday_reduced-basis_2002,quarteroni_reduced_2016}. In contrast to these so-called \emph{a posteriori} methods, \emph{a priori} approaches such as the Proper Generalised Decomposition (PGD) \parencite{chinesta_short_2011} build the reduced-order basis on the fly without the need for prior full-order computations. Initially introduced for space-time problems as the \emph{radial approximation} \parencite{Ladeveze85}, the PGD has since been applied to a broader range of parametric frameworks \parencite{ niroomandi_real-time_2013,neron_time-space_2015,lu_multi-parametric_2018,DabyHybrid}.

\bigbreak

Lately, we have also observed the successful use of machine learning methods (ML) to solve PDEs. The governing PDEs are usually incorporated into the loss function to ensure the model predictions are consistent with the known underlying physics. \AD{The physics-informed neural networks (PINNs) proposed in \parencite{raissi2019physics}, offer broader interpolation possibilities than the standard FEM.} \AD{Indeed,} in PINNs, classical fully connected neural networks are used to approximate a PDE solution. \AD{The training of such model is often} based on a loss function combining the \AD{influence of the physics (governing PDE residual, potential energy, etc.)}, boundary conditions, and optional measurement data, if they are available. The loss function is typically evaluated on a sufficiently large set of points randomly sampled within the computational domain.

\AD{However, a significant difficulty of PINNs is the imbalance among the loss terms when enforcing boundary conditions through penalization,} leading to unbalanced back-propagated gradients during the training \parencite{yu2018deepRitz,wang2021understanding}. Another approach to eliminating the problem with low accuracy in the boundary conditions consists in strongly imposing the boundary conditions. In \parencite{sukumar2022exact}, the neural network output is multiplied by approximate distance functions. The final solution is then obtained by adding the prescribed boundary values to the prediction that vanishes on the boundary. For the general domain, the computation of the distance function may be costly. So far, the feasibility of the proposed method has been demonstrated only on simple 2D domains. \AD{For tensor product meshes, Dirichlet boundary conditions can also be prescribed using a pre-computed function vanishing on the boundary of the domain as proposed in \parencite{omella2024r} where the physics is encoded by using the potential energy of the system as the loss function} \parencite{omella2024r}. \AD{Finally, the Dirichlet boundary condition can be accounted for by projection of the neural network interpolation onto a finite element mesh where the boundary conditions are prescribed and where the loss is evaluated. Such an approach is used in variational physics-informed neural networks (VPINNs) \parencite{berrone2022solving,berrone2022variational} where the loss contains the residuals of the variational form of the governing equations.}

\bigbreak

Recently, efforts have been made to combine the FEM and neural networks. Contrary to the original PINNs, the solution is typically obtained as a linear combination of predefined types of shape functions. The approximation capability of such models is restricted compared to neural networks; the solution is, however, produced in a more controlled manner. Due to the higher interpretability of the model parameters, the methods typically overcome some of the issues with the boundary conditions.

In the Neural-integrated meshfree (NIM) method proposed in \parencite{du2024neural}, the solution is interpolated using basis functions defined on the mesh-free discretization and nodal coefficient produced by a deep neural network. A fully connected neural network predicts the nodal coefficients. Dirichlet boundary conditions are incorporated through additional terms in the loss function.

In \parencite{badia2024finite}, the proposed method is built on the idea of VPINNS -- interpolating the NN solution to a space of Lagrangian shape functions and using a loss function based on a weak form of the PDE residuals. This time, the boundary conditions are strongly prescribed after the projection to the shape function space. A similar approach is used in the NEUFENET \parencite{khara2024neufenet}, which uses convolutional neural networks to predict nodal values associated with predefined shape functions. Dirichlet boundary conditions are not learned but are prescribed in a post-processing step. The neural network is trained using a loss function based on the potential energy of the problem. The method's applicability is illustrated on domains discretized by structured regular meshes. Due to the nature of the convolutional neural networks, the extension to general domains and unstructured meshes would not be straightforward. 

A slightly different approach is used in \AD{Embedded Finite element Neural Networks (EFENN) introduced with} the Hierarchical Deep learning Neural Networks (HiDeNN) framework \parencite{zhang2021hierarchical}. Contrary to the above-mentioned methods, where the shape functions are evaluated externally, in \AD{EFENN, the shape functions are embedded} directly within the neural network using a specific model architecture. 
This unified approach has several advantages. Firstly, the neural network's architecture is very sparse and is fully determined by the discretization of the problem. Second, all model parameters are interpretable. This property of the model can be exploited in various transfer learning strategies, for example, when training a model with a refined discretization network. The nodal coordinates and the nodal values are both parameters of the model, it is therefore possible to optimize them simultaneously. 
In \parencite{zhang2021hierarchical}, the HiDeNN model with Lagrangian shape functions and potential energy loss function was used on problems in 1D and 2D.

\bigbreak

\AD{Additionally, while some work relies on transfer learning and fine-tuning to limit the cost of the training stage for parametrized PDE \parencite{penwarden2021physics}, due to the black box-nature of the deep neural networks, a change of architecture means starting training from scratch. This is a major limitation of standard fully connected neural networks as the choice of the network architecture, the activation functions, and the optimization algorithm often remain arbitrary. Indeed, no general rule exists; therefore, finding the right combination for the specific problem remains a process of trial and error, and each change requires starting over and training from scratch the new model. However, we will see that the interpretability offered by the EFENN enables overcoming such difficulty.} 

\bigbreak

\AD{In this paper, we propose a new EFFEN framework named FENNI, which builds upon HiDeNN and introduces} \KS{several modifications and extensions to enhance its capabilities. \AD{We also further study factors affecting model accuracy and training efficiency.}} \AD{The FENNI framework} \KS{performs both the construction of shape functions and interpolation based on predicted nodal values \AD{and position, similar in purpose to the HiDeNN framework but with a distinct architecture}. \AD{Indeed, w}hile the architecture proposed in \parencite{zhang2021hierarchical} is elegant in 1D, its extension to higher dimensions presents challenges when considering unstructured meshes. 
	To address this, we propose an alternative architecture where shape functions are evaluated on the reference element. This approach enables us to define a unified architecture for higher spatial dimensions and varying orders of \AD{interpolation}. The reference-element-based architecture also allows us to leverage tabulated Gaussian quadrature points and weights for \AD{efficient and accurate integration during the} loss function evaluation.
	Finally, the concept of combined rh-adaptivity, proposed for the original HiDeNN model for 1D, is extended \AD{ and further applied to 2D, enabling localized mesh adaptivity. This newly introduced rh-adaptivity is hierarchical; the r-adaptivity sets a new Jacobian-based criterion that drives the h-adaptivity. Nodes are, therefore, added depending on the ongoing r-adaptivity using a red-green refinement strategy.}
}

\KS{
	In addition to the potential energy loss function used in the original HiDeNN method, which does not exist for mechanics problems with non conservative forces, we explore the use of two other loss functions: a residual loss typically used with PINNs and a loss function based on the weak formulation of the problem, both of which being defined for all mechanics problems, and studied here for the first time in the context of EFFEN. We examine how the choice of loss function influences the training process and solution accuracy on both fixed and adaptive meshes. Additionally, we explore the effect of different integration methods for evaluating the loss function. Finally, leveraging the model's interpretability, we introduce a multigrid training strategy and demonstrate its impact on the model robustness, accuracy and training time.
}

\AD{Part I of this paper, while not directly focusing on surrogate modelling aspects, paves the way to developments in versatile reduced-order modelling for designing patient-specific digital twins. This work aims to lay the elementary interpolation tools that will then be assembled into a high-dimension surrogate model in Part II \parencite{daby-seesaramFiniteElementNeural2024}.}

\AD{
	The outline of the paper is following. In Section \ref{sec:pb}, we introduce the mechanical problem used to illustrate the proposed FENNI framework as well as the different loss functions whose minimization leads to finding the solution to said problem. \cref{sec:method} presents the FENNI framework and its architecture dynamic optimization strategy during training.} \KS{The numerical experiments are presented in Section \ref{sec:res}. 
}

\section{Problem setting}\label{sec:pb}

\subsection{Definition of the problem} \label{sec:problem}

\begin{figure}[H]
	\centering
	\includegraphics[width=0.25\linewidth]{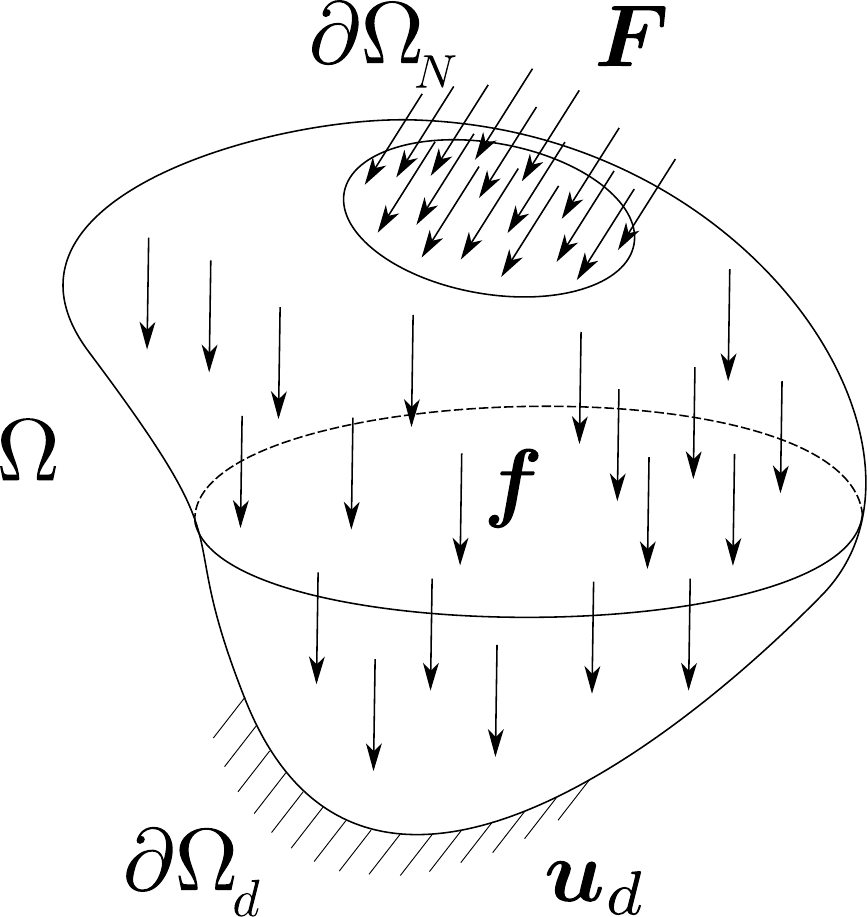}
	\caption{Reference problem}
	\label{fig:ref}
\end{figure}

\AD{To illustrate the proposed framework}, we consider a linear elastostatic problem in $\mathbb{R}^k$, where dimension $k=1,2,3$. 
Let us consider bounded connected domain $\Omega \in \mathbb{R}^k$ \AD{illustrated in \cref{fig:ref}}. The boundary of domain $\Omega$ consists of two subdomains $\partial \Omega = \partial\Omega_d \cup \partial\Omega_N$, such that $\partial\Omega_d \cap \partial\Omega_N = \emptyset$.

The governing equations read:
\begin{align}
	&\nabla \cdot \doubleunderline{\sigma}\left(\underline u\right) + \underline f = 0 &\quad \mathrm{in} \ \Omega, \label{eq:govern_1}\\
	& \underline{\underline{\sigma}}^T =  \underline{\underline{\sigma}} &\quad \mathrm{in} \ \Omega, \label{eq:govern_11}\\
	&\underline{n} \cdot \underline{\underline{\sigma}}\left(\underline u\right) = \underline t &\quad \mathrm{on} \ \partial\Omega_N, \label{eq:govern_2}\\
	&\underline{u} = \underline{u}_D &\quad \mathrm{on} \ \partial\Omega_d \label{eq:govern_3},
\end{align}
where $\underline u$ is displacement, $\doubleunderline \sigma$ is the Cauchy stress tensor, $\underline f$ is the prescribed body force, $\underline{u}_D$ and $\underline t$ are the displacement and surface force prescribed on the respective segments of the domain boundary. 
With the assumption of linear elasticity, the stress-strain relation is $\underline{\underline{\sigma}} = \underline{\underline{\underline{\underline{C}}}}: \underline{\underline{\varepsilon}}\left(\underline u\right)$, where $\doubleunderline{\varepsilon}\left(\underline u\right) := (\nabla \underline u + (\nabla \underline u)^T)/2$ and $\underline{\underline{\underline{\underline{C}}}}$ is the stiffness tensor.

\subsection{Formulation of loss functions} \label{sec:loss_functions}
Multiple loss functions can be designed such that their minimization will lead to \AD{the solution of the problem}. While all theoretically sound, the actual effectiveness of the loss functions, when used in a particular setting, can vary significantly.  
In this work, we use the potential energy loss function, a loss function containing the residuals of the governing equations typically used with classical PINNs and the loss function based on the weak formulation. 
The three types of loss function for governing equations (\ref{eq:govern_1}),(\ref{eq:govern_2}), (\ref{eq:govern_3}) are derived below.

\subsubsection{Potential energy loss function }
The potential energy loss, used in the Deep Ritz method \parencite{yu2018deepRitz}, was also used with the original HiDeNN model in \parencite{zhang2021hierarchical}.
The potential energy of \AD{the} system is defined as follows:
\begin{align}
	L_p = \frac{1}{2}\int_\Omega \doubleunderline{\sigma}\left(\underline u\right):\doubleunderline{\varepsilon}\left(\underline u\right) \mathrm{d} \Omega
	- \int_\Omega\underline f \cdot \underline u \mathrm{d} \Omega
	- \int_{\partial\Omega_N} \underline t \cdot \underline u \mathrm{d} A,
	\label{eq:loss_general_potential}
\end{align}
where the displacement $\underline{u}$ is the output of the model.
We remark that, contrary to most typically used loss functions, $L_p \neq 0$ at the minimum. This property needs to be taken into account when designing the convergence criterion. \AD{Moreover it must be noted that such a loss is limited to conservative loading, for which a potential energy can be defined. } 

\subsubsection{Residual loss function}
Alternatively, we can use the residual loss function proposed for the original PINNs \parencite{raissi2019physics,wang2021understanding}. With our model, the loss is defined as a combination of the residuals in equations (\ref{eq:govern_1}) and (\ref{eq:govern_2}). The residuum of equation (\ref{eq:govern_3}) is not included in the loss function because the Dirichlet boundary conditions for $\underline u$ -- the output of the model -- can be strongly enforced within the model, as will be described further. \AD{Similarly, the symmetry of the stress tensor is strongly enforced.  }
For a set of training points $X_{\Omega} = \{x_i\ |\ i=0,\dots,n_{\Omega}\}$ sampled inside domain $\Omega$, and a set of points $X_{\partial\Omega} = \{x_i\ |\ i=0,\dots,n_{\partial\Omega}\}$ sampled on the boundary $\partial\Omega_N$, the loss is defined as a combination of two loss terms:
\begin{align}
	\tilde L_r = \frac{\lambda_1}{n_{\Omega}}\sum_{i=1}^{n_{\Omega}}\left( \nabla \cdot \doubleunderline{\sigma}\left(\underline u\right) - \underline f \right)^2 
	+ \frac{\lambda_2}{n_{\partial\Omega}}\sum_{i=1}^{n_{\partial\Omega}}\left( \underline{n} \cdot \underline{\underline{\sigma}}(u) - \underline t \right)^2
\end{align}

To benefit from the strongly prescribed boundary conditions also for the stress, we can reformulate the problem. Instead of computing the stress based on the predicted displacement $\underline u$
as $\underline{\underline{\sigma}}\left(\underline u\right) = \underline{\underline{\underline{\underline{C}}}}: \doubleunderline{\varepsilon}\left(\underline u\right)$, we \AD{turn to a mixed formulation and} train two models simultaneously -- one for $\underline u$ and one for $\hat{\underline{\underline{\sigma}}}$. This way, stress and displacement are independent variables, and both boundary conditions given by equation (\ref{eq:govern_2}) and (\ref{eq:govern_3}) can be prescribed strongly within the corresponding models. The loss function then contains one residuum term and a compatibility term minimizing the difference between the stress predicted by the model $\hat{\doubleunderline{\sigma}}$ and stress computed \AD{from} the displacement $\doubleunderline{\sigma}\left(\underline u\right)$:
\begin{align}
	L_r = \frac{\lambda_1}{n_{\Omega}}\sum_{i=1}^{n_{\Omega}}\left( \nabla \cdot \hat{\doubleunderline{\sigma}} - \underline f \right)^2 
	+ \frac{\lambda_1}{n_{\Omega}}\sum_{i=1}^{n_{\Omega}}\left( \doubleunderline{\sigma}\left(\underline u\right) - \hat{\doubleunderline{\sigma}} \right)^2 
	\label{eq:loss_general_resid}
\end{align}
Another advantage of loss function $L_r$ compared to $\tilde L_r$ is that it is evaluated on only one set of training points. We, therefore, avoid the problem of balancing the number of points samples in domain $\Omega$ and on the boundary $\partial\Omega_N$.

\subsubsection{Weak formulation loss function}
In the spirit of the FEM, the loss function can also be formulated using the weak formulation of the governing equations. For example, this type of loss function was used in the VPINN \parencite{berrone2022solving}. The advantage of using this type of loss function \AD{within our framework}, compared to the deep neural networks in VPINN, is that the model loss can be evaluated without applying any projection to the model output. By doing so, we also remove a potential source of error introduced by the projection step.

To derive the loss function, we start by deriving the weak formulation of equation (\ref{eq:govern_1}), 
\begin{align}
	W:=\underbrace{\int_\Omega \doubleunderline{\sigma}\left(\underline u\right):\doubleunderline{\varepsilon}\left(\underline u^*\right) \mathrm{d}\Omega}_
	{W_{\mathrm{int}}(\underline u;\underline u^*)}
	- \underbrace{ \left( \int_{\partial\Omega_N} \underline t \cdot \underline u^* \mathrm{d}A+ \int_{\Omega} \underline f \cdot \underline u^* \mathrm{d}\Omega \right) }_
	{W_{\mathrm{ext}}(\underline u;\underline u^*)}
	= 0,
	\quad \forall u^* \in H_1^0(\Omega).
	\label{eq:weak_form}
\end{align}

The problem is solved on finite-dimension space defined by a set of basic functions $\mathcal{B} = \{u^*_i\}_{i=0}^n$ leading to 
\begin{align}
	W\left(\underline u, \underline u^*\right) = 0, \quad \forall \underline u^* \in H_{1,n}^0,
	\label{eq:weak_resid_2}
\end{align}
where $H_{1,n}^0 =\{ \underline v \in \mathrm{Span}\left(\mathcal{B}\right)\ |\ \underline v = 0 \mathrm{\ on\ } \partial \Omega_D \}$. 

The loss function is constructed based on the following implication:
\begin{align}
	\sum_{i=1}^{n} \left(W\left(\underline u,\underline u_i^*\right)\right)^2 = 0, \Rightarrow W\left(\underline u,\underline u^*\right) = 0, \quad \forall \underline u^* \in H_{1,n}^0.
	\quad
\end{align}

Therefore, by minimizing the following loss function

\begin{align}
	L_w =  \sum_{i=1}^{n} \left(  
	\int_\Omega \doubleunderline{\sigma}\left(\underline u\right):\doubleunderline{\varepsilon}\left(\underline u_i^*\right) \mathrm{d}\Omega - \int_{\partial\Omega_N} \underline t \cdot \underline u_i^* \mathrm{d}A - \int_\Omega\underline f \cdot \underline u_i^* \mathrm{d}\Omega
	\right)^2,
	\label{eq:loss_general_weak}
\end{align}
we ensure that equation (\ref{eq:weak_resid_2}) is satisfied. \AD{Note that while this formulation already yields good results, a continuous formulation of such weak equilibrium gap with a consistent discretization has been proposed by} \cite{genetFiniteStrainFormulation2023a} \AD{and used by} \cite{peyrautFiniteStrainFormulation2024} \AD{for parameters identification and could be used here as well.}

\section{Methods}\label{sec:method}

\AD{In this section, we first recall the initial EFENN framework introduced with the HiDeNN by \cite{zhang2022hidenn} before exposing the improvements brought by the proposed FENNI framework	}

\subsection{Recall of the existing EFENN - HiDeNN framework (interpolation-layer-based formulation)}

\KS{We start by recalling the construction of shape functions in the original 1D HiDeNN model \parencite{zhang2022hidenn}. However, we introduce a new formalism relying on Assembly Layer, which enables a unified handling of different orders of shape functions.}

Let us consider a 1D domain $\Omega$ discretized by a mesh of $np$ nodes and $ne = np-1$ elements. 
Similarly to \parencite{zhang2021hierarchical}, the first step in the construction of shape functions is the linear building block:
\begin{align}
	\mathrm{Linear\ Block}\left(x, x_a, x_b, y_a, y_b\right) = \left(y_b-y_a\right)\mathrm{ReLU}\left(-\frac{1}{x_b-x_a}\mathrm{ReLu}\left(-x+x_b\right)+1\right) + y_a,
\end{align}
where $\mathrm{ReLu}(x) = max\{0,x\}$ is the rectified linear unit activation function.
The linear building block generates piece-wise linear functions that are subsequently used to construct the shape functions. 
For each element $e_i$, $i=0,\dots,ne-1$, the two piece wise linear auxiliary functions functions $\tilde N_i^0, \tilde N_i^1$ are constructed:
\begin{align}
	\tilde N_i^0 &= \mathrm{Linear\ Block}\left(x, x_i, x_{i+1}, 1, 0\right) \\
	\tilde N_i^1 &= \mathrm{Linear\ Block}\left(x, x_i, x_{i+1}, 0, 1\right).
\end{align}
The function responsible for the generation of the auxiliary function for a given element is named Element Block:
\begin{align}
	\mathrm{Element\ Block}\left(x,x_i,x_i+1\right)=\left(\tilde N_i^0,\tilde N_i^1\right).
\end{align}

\begin{figure}[h!]
	\centering
	\begin{subfigure}{0.75\textwidth}
		\includegraphics[width=\textwidth]{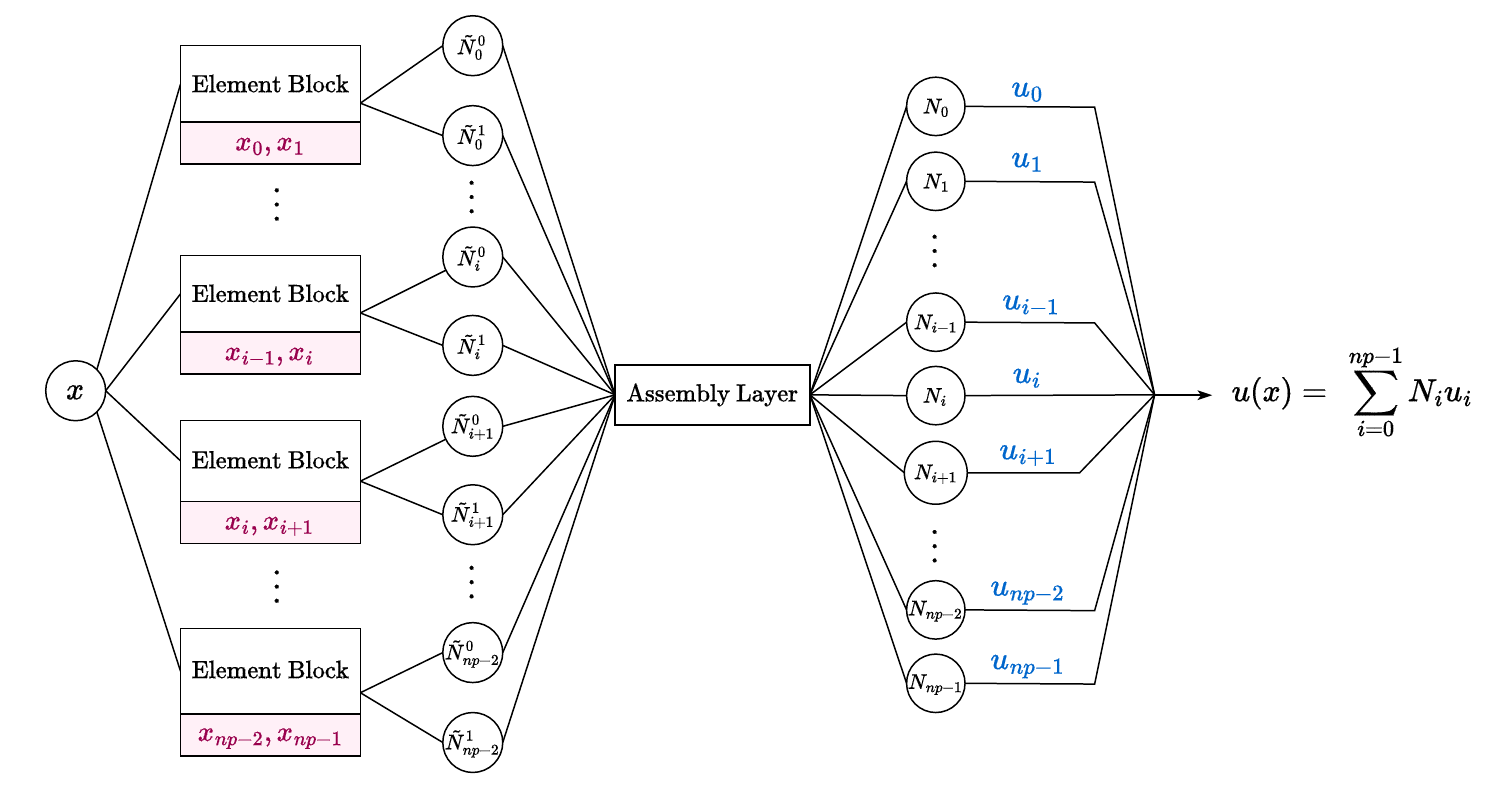}
		\caption{Interpolation-layer-based implementation}
		\label{fig:impl_interpol}
	\end{subfigure}
	\par\bigskip
	\begin{subfigure}{0.75\textwidth}
		\includegraphics[width=\textwidth]{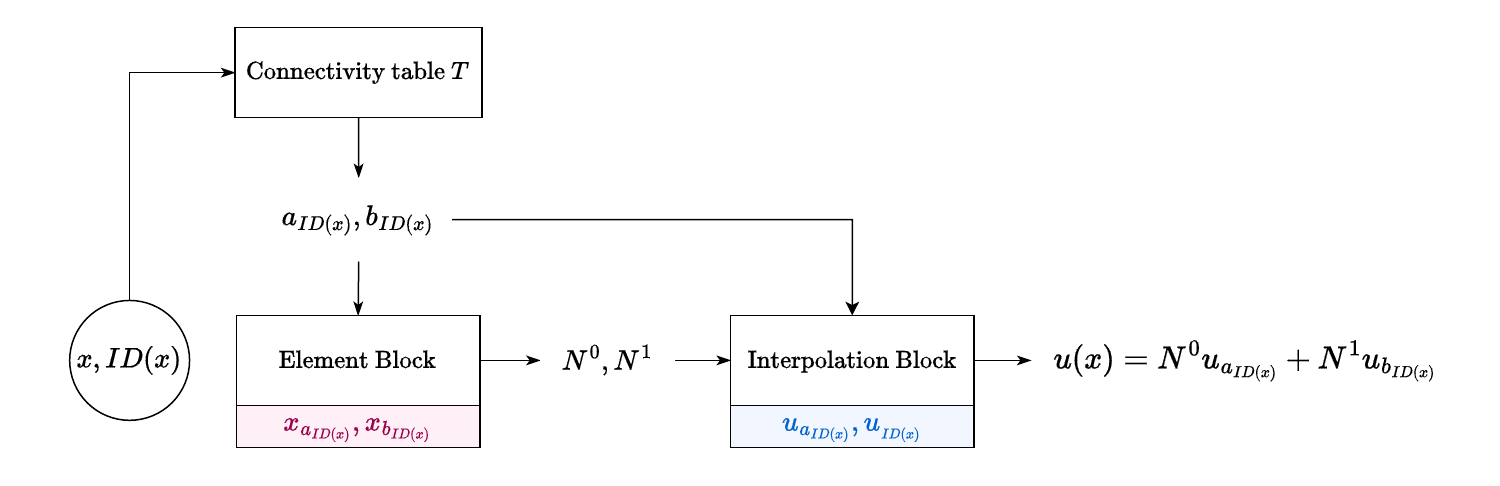}
		\caption{Reference-element-based implementation}
		\label{fig:impl_ref}
	\end{subfigure}
	\caption{Two versions of the model implementation, defined for domain $\Omega$ discretized by a mesh with $np$ nodes and linear shape functions. In the interpolation-layer-based implementation (a), all shape functions are assembled and evaluated at the input point $x$. In the reference-element-based implementation (b), only the shape functions with element $e_{\text{ID}(x)}$ in their support are evaluated. If the problem is solved on fixed mesh, the nodal values $U_{np}$ are trainable, and nodal coordinates $X_{np}$ are fixed. For r-adaptivity, the nodal coordinates $X_{np}$ are also made trainable, which leads to simultaneous optimization of nodal values and their position. }
	\label{fig:implementation}
\end{figure}

The shape functions $N_i,\ i = 0,\dots,np-1$, are obtained by assembly of the auxiliary functions. The non-zero tails of functions $\tilde N_i^0, \tilde N_i^1$ are eliminated in the assembly process. 
For linear shape functions, each shape function associated with the inner node is constructed by assembling two auxiliary functions. Only one auxiliary function is needed for the shape functions of the boundary nodes. The assembly procedure is following:
\begin{align}
	N_{i} &= \tilde N_i^0 + \tilde N_{i-1}^1 -1,   \quad i=1,\dots,np-2,\\
	N_0 &= \tilde N_0^0,\\
	N_{np-1} &= \tilde N_{np-2}^1
\end{align}
The operations can be written in a compact form as $\underline N = \mathbb{A}\tilde{\underline N} + \underline b$, where $N \in \mathbb{R}^{np}$ is the vector consisting of the shape functions and $\tilde{\underline N} \in \mathbb{R}^{2\left(np-1\right)}$ is the vector of the auxiliary functions, obtained by concatenating outputs of all ElementBlocks. \KS{Because the assembly operation follows the form $\mathbb{A}{\underline x} + \underline b$, resembling a classic neural network layer, we refer to this block as the Assembly Layer. Its position within the model architecture is shown in the diagram in Figure \ref{fig:impl_interpol}.}

\begin{figure}[h]
	\centering
	\begin{subfigure}{0.35\textwidth}
		\includegraphics[width=\textwidth]{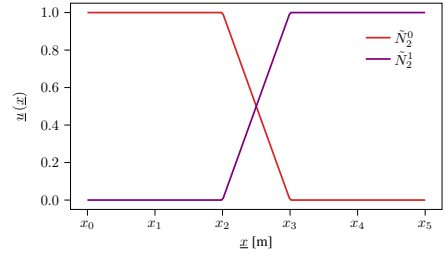}
		\caption{The two auxiliary functions defined for one element.}
	\end{subfigure}

	\hspace{1cm}
	\begin{subfigure}{0.4\textwidth}
		\includegraphics[width=\textwidth]{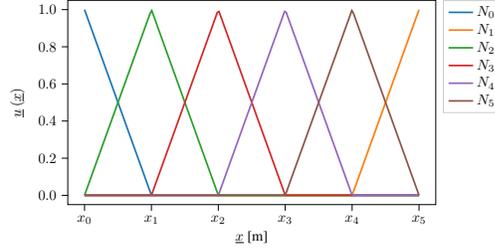}
		\caption{Assembled global shape functions.}
	\end{subfigure}
	\caption{The construction of linear shape functions. 
		The global linear shape functions (b) are assembled by combining the auxiliary function for all elements. Two auxiliary functions are defined for each element (a). All shape functions associated with interior nodes of the domain are constructed by combining the auxiliary shape functions from two neighboring elements and eliminating the non-zero tails.}
	\label{fig:interpolation_lin}
\end{figure}

The two types of auxiliary functions defined for each element and the assembled linear shape functions are shown in Figure \ref{fig:interpolation_lin}.

\begin{figure}[h]
	\centering
	\begin{subfigure}{0.355\textwidth}
		\includegraphics[width=\textwidth]{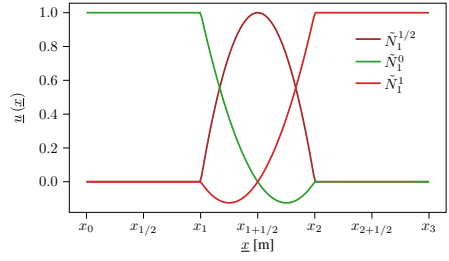}
		\caption{The three auxiliary functions defined for one element.}
	\end{subfigure}

	\hspace{1cm}
	\begin{subfigure}{0.425\textwidth}
		\includegraphics[width=\textwidth]{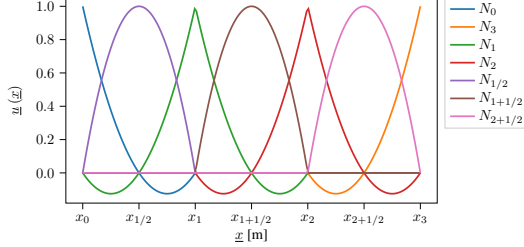}
		\caption{Assembled global shape functions.}
	\end{subfigure}
	\caption{The global quadratic shape functions (b) are assembled by combining the auxiliary function for all elements. Three auxiliary functions are defined for each element (a). Each shape function associated with an interior node is constructed by combining two auxiliary shape functions and removing the non-zero tails. Only one auxiliary function is used for the shape function associated with the element mid-point.}
	\label{fig:interpolation_quad}
\end{figure}

For quadratic shape functions, we define three auxiliary functions $\tilde N_i^0, \tilde N_i^1, \tilde N_i^{1/2}$ for each element $e_i$. The auxiliary functions are again defined inside the Element block using the linear building block:
\begin{align}
	\tilde N_i^0 &= \frac{\mathrm{Linear\ Block}(x, x_i, x_{i+1}, 1, x_{i+1} - x_i)
		\cdot \mathrm{Linear\ Block}\left(x, x_i, x_{i+1}, x_{i+1} - x_{i}, 0\right)}
	{(x_{i+1} - x_{i+1})\left(x_{i+1} - x_{i+1/2}\right)},\\
	\tilde N_i^1 &= \frac{\mathrm{Linear\ Block}(x, x_i, x_{i+1},x_{i+1/2} - x_{i} , x_{i+1/2} - x_{i+1})\cdot \mathrm{Linear\ Block}(x, x_i, x_{i+1}, x_{i+1} - x_{i}, 0)}
	{(x_{i} - x_{i+1/2})\left(x_{i} - x_{i+1}\right)},\\
	\tilde N_i^{1/2} &= \frac{\mathrm{Linear\ Block}\left(x, x_i, x_{i+1},0 , x_{i+1} - x_{i}\right)\cdot 
		\mathrm{Linear\ Block}\left(x, x_i, x_{i+1}, x_{i+1} - x_{i}, 0\right)}
	{(x_{i+1/2} - x_{i})\left(x_{i+1/2} - x_{i+1}\right)}.
\end{align} 
Node $x_{i+1/2}$ is positioned in the middle of element $e_i$, that is $x_{i+1/2} = (x_i + x_{i+1})/2$.
The shape functions are obtained by assembly of these auxiliary functions. Contrary to the linear shape functions, the number of auxiliary functions used to construct a shape function is different. Furthermore, only two of the three auxiliary functions $\tilde N_i^0, \tilde N_i^1, \tilde N_i^{1/2}$ have non-zero tails that need to be eliminated during the assembly.   
The assembly procedure is, therefore following:
\begin{align}
	N_i &= \tilde N_{i-1}^1 + \tilde N_{i}^0 - 1, \quad i=1,\dots,np-2\\
	N_{i+1/2} &= \tilde N_i^{1/2}, \quad i=0,\dots,np-2,\\
	N_0 & = \tilde N_{0}^0,\\
	N_{np-1} & = \tilde N_{np-2}^1.
\end{align}
The three types of auxiliary functions defined for each element and the assembled quadratic shape functions are shown in Figure \ref{fig:interpolation_quad}.
\KS{The higher-order shape function would \AD{be} generated \AD{similarly}, making the Assembly Layer the only component that requires specific adaptation for each order of shape functions. In 1D, where the elimination of non-zero tails is straightforward due to a fixed number of neighboring elements, the general form $\mathbb{A}\underline x + \underline b$ is maintained. The same approach, however, cannot be easily extended to general meshes in 2D.}

To explain the interpolation process, we consider linear shape functions for simplicity. When the model is evaluated at point $x$, the output of the Assembly Layer is a vector of dimension $np$ consisting of the shape functions evaluated at point $x$: $\left(N_0(x),N_1(x),\dots,N_{np-1}(x)\right)$.
Therefore, the interpolation can be performed using a single linear layer with $np$ input nodes and one output node, zero bias, and the weights corresponding to the nodal values $u_i,\ i=0,\dots,np-1$. The output of the model for input $x$ is then equal to $\sum_{i=0}^{np-1} N_i(x)u_i$.

The evaluation of integral loss function in a form:
\begin{align}
	L\left(\underline{\theta}\right) = \int_\Omega l\left(x,\underline{\theta}\right) \ \mathrm{d}x,
\end{align}
where $\underline{\theta}$ represents the trainable parameters of the model, is done using the trapezoidal rule.
With notation $\Omega = [a,b]$, the integral loss is approximated by:
\begin{align}
	\int_a^b l\left(x,\underline{\theta}\right) \ \mathrm{d}x \approx \sum_{k=1}^{N_{int}}\frac{l(x_{k-1},\underline{\theta})+l\left(x_{k},\underline{\theta}\right)}{2}\left({x_{k}-x_{k-1}}\right),
\end{align}
where $a=x_1 < x_2 < \dots < x_{N_{int}} = b$ are the sampling points.
To compute the loss, the model needs to be evaluated in all of the sampling points.

\KS{The simplicity of the interpolation layer is the key advantage of this model implementation. However, several of its limitations need to be pointed out. Firstly, this architecture relies on the feasibility of assembling the global shape functions. As mentioned above, considering the nature of the auxiliary function generated by the linear building block, eliminating the non-zero tails is nontrivial in higher dimensions.
	Additionally, although the final sum $\sum_{i=0}^{np-1}u_i N_i(x)$ has very few non-zero terms, this implementation requires each input point $x$ to be processed in all $np-1$ element blocks and subsequently to assemble and evaluate all the shape functions\AD{, which suggests that the strategy could be made more efficient}. }

\subsection{The Finite Element Neural Network Interpolation (FENNI) framework} \label{sec:model}

The architecture of the proposed model is entirely determined by the number of nodes and elements used to discretize the domain and the order of shape functions used for the interpolation. 
Let us consider a domain $\Omega$ discretized by a mesh with $np$ degrees of freedom.
The model based on this discretization will be denoted by $M_{np}$ and it operates with two sets of parameters: the nodal values $U_{np} = \{u_i \ |\ i = 0, \dots, np-1\}$ and the nodal coordinates $X_{np} = \{x_i \ |\ i = 0, \dots, np-1\}$. 

\begin{figure}[h]
	\centering
	\begin{subfigure}{0.32\textwidth}
		\includegraphics[width=\textwidth]{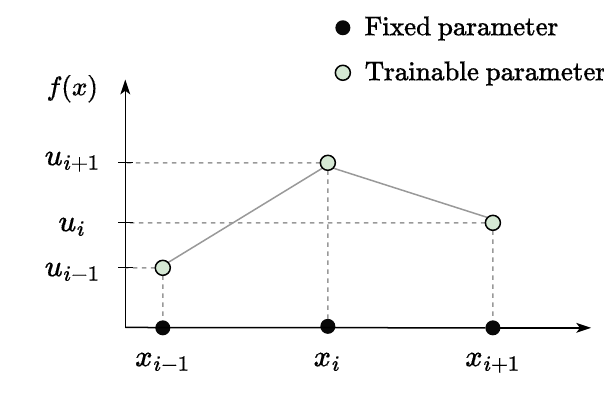}
		\caption{Fixed mesh}
	\end{subfigure}
	\begin{subfigure}{0.32\textwidth}
		\includegraphics[width=\textwidth]{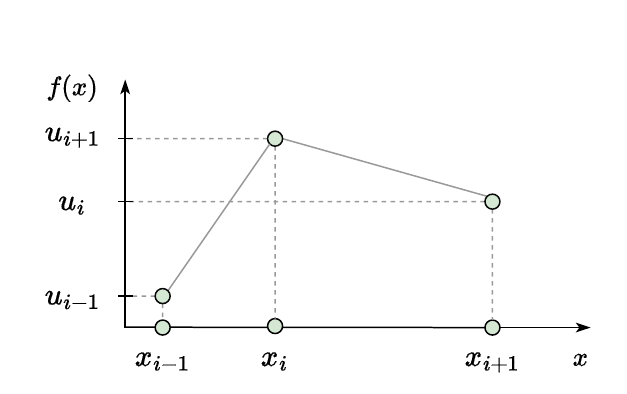}
		\caption{r-adaptivity}
	\end{subfigure}
	\begin{subfigure}{0.32\textwidth}
		\includegraphics[width=\textwidth]{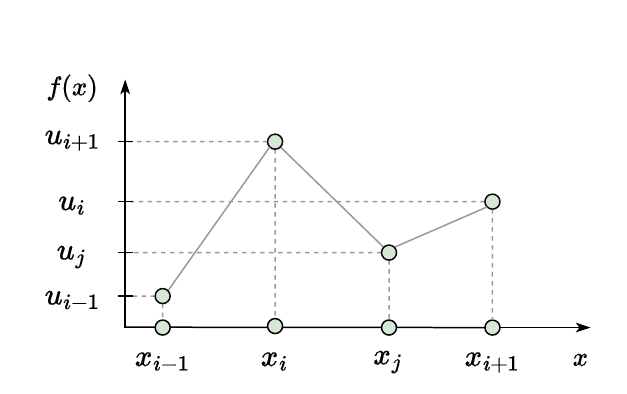}
		\caption{rh-adaptivity}
	\end{subfigure}
	\caption{The effect of trainable nodal values and nodal coordinates in 1D. }
	\label{fig:1D_basic}
\end{figure}

During the training, the model parameters $U_{np},\ X_{np}$ are optimized to minimize the loss function. The architecture of the model allows us to solve the problem in several settings. 
If only the nodal values $U_{np}$ are trainable, the training corresponds to solving the problem on a fixed mesh. Alternatively, if both the nodal values and coordinates are trainable, the training process consists of simultaneous optimization of nodal values and the position of the nodes. Therefore, by making $U_{np}$ and $X_{np}$ trainable, we achieve r-adaptivity.
Finally, adding new parameters during the training corresponds to adding new nodes in the mesh. By allowing the addition of new parameters, \emph{i.e.,} extension of sets $U_{np}$ and $X_{np}$, we achieve h-adaptivity. Furthermore, the rule based on which the nodes are added to specific regions can be linked to the movement of the nodes, resulting in combined rh-adaptivity. The three training options are illustrated in Figure \ref{fig:1D_basic}.

The proposed model can be viewed as a sparse neural network where values of all weights and biases are interpretable. This interpretability can be leveraged, for example, in prescribing Dirichlet boundary conditions. The nodal values corresponding to the boundary nodes can be fixed on the values prescribed by boundary conditions and are not modified during the training.

\KS{Finally, the interpretability of the parameters enables knowledge transfer between models trained on different discretizations of the same domain. As demonstrated later, this approach allows a significant portion of the training to be performed on a coarser mesh, resulting in improved training efficiency.}
\KS{Alternatively, instead of relying on a previously trained model, a new model can be initialized using solutions generated by classical FEM solvers.}
This is not the case for classical architectures of physically informed neural networks, where the transfer learning strategies typically rely on freezing the weights and biases of several pre-trained layers \parencite{xu2023transfer}.

\subsubsection{Reference-element-based formulation in 1D}

In this novel implementation, we address some of the limitations of the interpolation-layer-based implementation, mainly the extendability from 1D to general mesh in 2D.
In the reference-element-based implementation, we drop the concept of global shape functions and a single layer performing the interpolation. For an input point $x$, the shape functions are defined only on the element where $x$ is located. Thus, we can avoid the problem of eliminating any non-zero contributions coming from neighboring elements.
The same approach is used for any order of shape functions; therefore, we focus only on linear shape functions.

The input of the model is a tuple $\left(x, \text{ID}(x)\right)$ - a coordinate of the point and ID of the element in which the point is localized.
We again consider a domain $\Omega$ discretized by a mesh with $np$ nodes and $ne$ elements. In addition to the two sets of potentially trainable parameters -- nodal values $U_{np}$ and nodal coordinates $X_{np}$ -- the model also uses a connectivity table $T = \left\{\left[a_i, b_i\right] \ |\ i = 0,\dots,np-2\right\}$, where $a_i, b_i$ are identifiers of nodes in element $e_i$. \emph{I.e.}, $x_{a_i}$ and $x_{b_i}$ are the coordinates of nodes in element $e_i$ and the nodal values at these points are $u_{a_i}$ and $u_{b_i}$. The connectivity table is built based on the discretization mesh during model creation. 

\begin{figure}[h]
	\centering
	\begin{subfigure}{0.35\textwidth}
		\includegraphics[width=\textwidth]{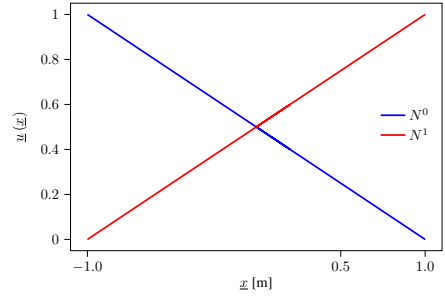}
		\caption{Shape functions defined on the reference element.}
	\end{subfigure}

	\hspace{1cm}
	\begin{subfigure}{0.428\textwidth}
		\includegraphics[width=\textwidth]{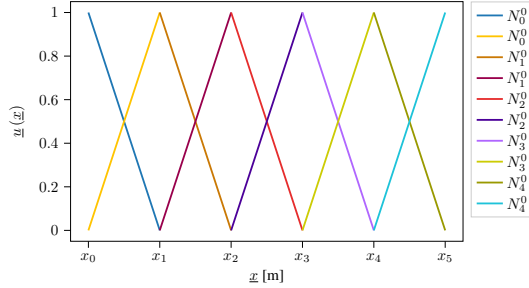}
		\caption{Shape functions restricted to one element.}
	\end{subfigure}
	\caption{Contrary to the interpolation-layer-based implementation, the shape functions are defined on the reference element (a), and the global shape functions are not constructed. For the interpolation, shape functions restricted to one element are used (b).}
	\label{fig:1D_ref_elem}
\end{figure}

The reference element in 1D is interval $e = [-1,1]$. For each input tuple $(x, \text{ID}(x))$, the position of $x$ in the reference element $e$ are computed as:
\begin{align}
	x^{\text{ref}} = 2 \frac{x-x_{a_{\text{ID}(x)}}}{x_{b_{\text{ID}(x)}}-x_{a_{\text{ID}(x)}}}-1,
\end{align}
where $x_{a_{\text{ID}(x)}}, x_{b_{\text{ID}(x)}}$ are the boundary nodes of element $e_{\text{ID}(x)}$. The coordinates of boundary nodes are identified based on the table $T$.
The shape functions are evaluated based on the reference coordinates $x^{\text{ref}}$:
\begin{align}
	N^0 = -0.5 \cdot x^{\text{ref}} + 0.5,\\
	N^1 = 0.5 \cdot x^{\text{ref}} + 0.5.
\end{align}
To get the interpolated value $u(x)$, the shape functions $N_{\text{ID}(x)}^0, N_{\text{ID}(x)}^1$ are combined with the corresponding nodal values found based on the connectivity table:
\begin{align}
	u(x) = N^0\left(x^{\text{ref}}\right) \cdot u_{a_{\text{ID}(x)}} + N^1\left(x^{\text{ref}}\right) \cdot u_{b_{\text{ID}(x)}} 
\end{align}
The construction of the interpolation based on the shape functions restricted to one element is illustrated in Figure \ref{fig:1D_ref_elem}.

Contrary to interpolation-layer-based implementation, the global shape functions are never fully assembled, and we operate only with their restrictions to a given element. The implementation is based on providing the element ID as an input parameter alongside the coordinates of point $x$. Thanks to the known localization of $x$, only the necessary shape functions are constructed and evaluated. 

The integral loss function:
\begin{align}
	L\left(\underline{\theta}\right) = \int_\Omega l\left(x,\underline{\theta}\right) \ \mathrm{d}x,
\end{align}
was evaluated using the Gauss quadrature rule.
With the $N_g$-point Gauss quadrature, the integral loss is approximated by:
\begin{align}
	\nonumber
	\int_a^b l\left(x,\underline{\theta}\right) \ \mathrm{d}x &= \sum_{j=1}^{ne} \int_{x_{a_j}}^{x_{b_j}} l\left(x,\underline{\theta}\right) \ \mathrm{d}x
	= \sum_{j=1}^{ne} \int_{-1}^{1} l\left(x\left(x^{\text{ref}}\right),\underline{\theta}\right) \ \frac{\mathrm{d}x}{\mathrm{d}x^{\text{ref}}}\mathrm{d}x^{\text{ref}} \\
	&\approx \sum_{j=1}^{ne} \sum_{i=1}^{N_g} w_i l\left( \frac{x_{b_j} - x_{a_j}}{2} x^g_i + \frac{x_{a_j} + x_{b_j}}{2}, \underline{\theta} \right) \frac{x_{b_j} - x_{a_j}}{2},
\end{align}
where $x^g_i$, $w_i$, $i=0,\dots,N_g-1$ are the Gauss quadrature points and weights for interval $[-1,1]$.
The model needs to be evaluated at the predefined quadrature points in each element to compute the loss.

\subsubsection{Extension of the reference-element-based implementation to 2D}

One of the key advantages of the reference-element-based implementations is its straightforward extension to higher dimensions. In this section, we describe the extension to 2D, but the situation in 3D would be analogical.

Analogically to 1D, the input of the model consists of coordinate $\underline x$ and ID of the element in which $\underline x$ is localized -- a tuple $\left(\underline x, \text{ID}\left(\underline x\right)\right)$.
We consider a domain $\Omega$ discretized by a triangular mesh with $np$ nodes and $ne$ elements. The model operates with the two sets of potentially trainable of parameters: nodal values $U_np = \{\underline u_i \ |\ i=0,\dots,np-1\}$ and nodal coordinates $X_np = \{\underline x_i\ |\ i=0,\dots,np-1\}$. Additionally, we define the connectivity table $T = \{[a_i, b_i, c_i] \ |\ i = 0,\dots,ne-1 \}$, where $a_i, b_i, c_i$ are identifiers of nodes in element $e_i$. 
In 2D, the reference element is a triangle $e = \{ \underline y=(y_1,y_2)\ |\ y_1\geq 0,\ y_2\geq 0,\ y_1+y_2 \leq 1 \}$.
For each input tuple $(\underline x, \text{ID}\left(\underline x\right))$, the position of point $x$ in the reference element are computed as:
\begin{align}
	\begin{pmatrix}
		x^{\text{ref}}\\
		y^{\text{ref}}\\
		1 - x^{\text{ref}} - y^{\text{ref}}
	\end{pmatrix} = 
	\begin{pmatrix}
		x_{a_{\text{ID}(x))}} & x_{b_{\text{ID}(x))}} & x_{b_{\text{ID}(x))}}\\
		y_{a_{\text{ID}(x))}} & y_{b_{\text{ID}(x))}} & y_{b_{\text{ID}(x))}}\\
		1 & 1 & 1
	\end{pmatrix}^{-1}
	\begin{pmatrix}
		x\\
		y\\
		1
	\end{pmatrix} 
\end{align}
where ${\underline{x}}_{a_{\text{ID}(x))}} = \left(x_{a_{\text{ID}(x))}},y_{a_{\text{ID}(x))}}\right),\ {\underline{x}}_{b_{\text{ID}(x))}} = \left(x_{b_{\text{ID}(x))}},y_{b_{\text{ID}(x))}}\right),\ {\underline{x}}_{c_{\text{ID}(x))}} = \left(x_{c_{\text{ID}(x))}},y_{c_{\text{ID}(x))}}\right)$ are the coordinates of the nodes of element $e_{\text{ID}\left(\underline x\right)}$. 
The shape functions restricted to element $e_{\text{ID}(x)}$ are computed based on the reference coordinates $\underline{x}^{\text{ref}}$:
\begin{align}
	N_{\text{ID}(x)}^0 &= x^{\text{ref}},\\
	N_{\text{ID}(x)}^1 &= y^{\text{ref}},\\
	N_{\text{ID}(x)}^2 &= 1- x^{\text{ref}} - y^{\text{ref}}.
\end{align}
Analogically to the 1D case, the interpolated value of $\underline u(\underline x)$ is obtained using the corresponding nodal values:
\begin{align}
	\underline u(\underline x) = N_{\text{ID}(x)}^0 \cdot \underline u_{a_{\text{ID}(x)}} + N_{\text{ID}(x)}^1 \cdot \underline u_{b_{\text{ID}(x)}} + N_{\text{ID}(x)}^2 \cdot \underline u_{c_{\text{ID}(x)}}. 
\end{align}
The mapping between the elements and the reference element is always assumed to be linear. Therefore, the same approach can be used for any order of shape functions.
Analogically to the 1D case, the r-adaptivity is achieved by making the nodal coordinates trainable parameters. The coordinates of boundary nodes of domain $\Omega$ are kept fixed to preserve the shape of the domain. 

We consider a loss function in the form of an integral over domain $\Omega$:
\begin{align}
	L\left(\underline{\theta}\right) = \int_\Omega l\left(\underline x,\underline{\theta}\right) \ \mathrm{d}\underline x.
\end{align}
Similarly to 1D, we can benefit from the reference-element-based implementation and using the tabulated values of quadrature points and weights.
With the $N_g$-point Gauss quadrature, the integral loss is approximated as follows:
\begin{align}
	\nonumber
	\int_\Omega l\left(\underline x,\underline{\theta}\right) \ \mathrm{d}\underline x =
	\sum_{j=1}^{ne} \int_{e_j} l\left(\underline x,\underline{\theta}\right) \ \mathrm{d}\underline x &=
	\sum_{j=1}^{ne} \int_{e_{\text{ref}}} l\left(\underline x\left(\underline x^{\text{ref}}\right),\underline{\theta}\right) \left| J_j\right| \mathrm{d}\underline{x}^{\text{ref}} \\
	&\approx \sum_{j=1}^{ne} \sum_{i=1}^{N_g} w_i l\left( \underline x\left(\underline x^g_i,e_j\right) , \underline \theta \right) \left| J_j\right|
	\label{eq:int_loss_ref_elem}
\end{align}
For element $e_j$ with nodal coordinates $\underline x_{a_j}$, $\underline x_{b_j}$, $\underline x_{c_j}$, the physical coordinates of $i$-th quadrature point are:
\begin{align}
	\underline x(\underline x^g_i,e_j) = (x^g_ix_{a_j} + y^g_ix_{b_j} + (1-x^g_i-y^g_i)x_{c_j}, y^g_iy_{a_j} + y^g_iy_{b_j} + (1-x^g_i-y^g_i)y_{c_j})
\end{align}
and the Jacobian determinant $J_j$
\begin{align}
	J_j= \left(x_{a_j} - x_{c_j}\right)\left(y_{b_j} - y_{c_j}\right) - \left(x_{b_j} - x_{c_j}\right)\left(y_{a_j} -y_{c_j}\right)
\end{align}
is double the area of element $e_j$; \AD{a propriety that will be used to set the h-adaptivity criterion in} \cref{sec:mesh_adaptivity}.

\subsubsection{Mesh adaptation strategy}\label{sec:mesh_adaptivity}
\begin{figure}[h!]
	\centering
	\includegraphics[width=0.8\textwidth]{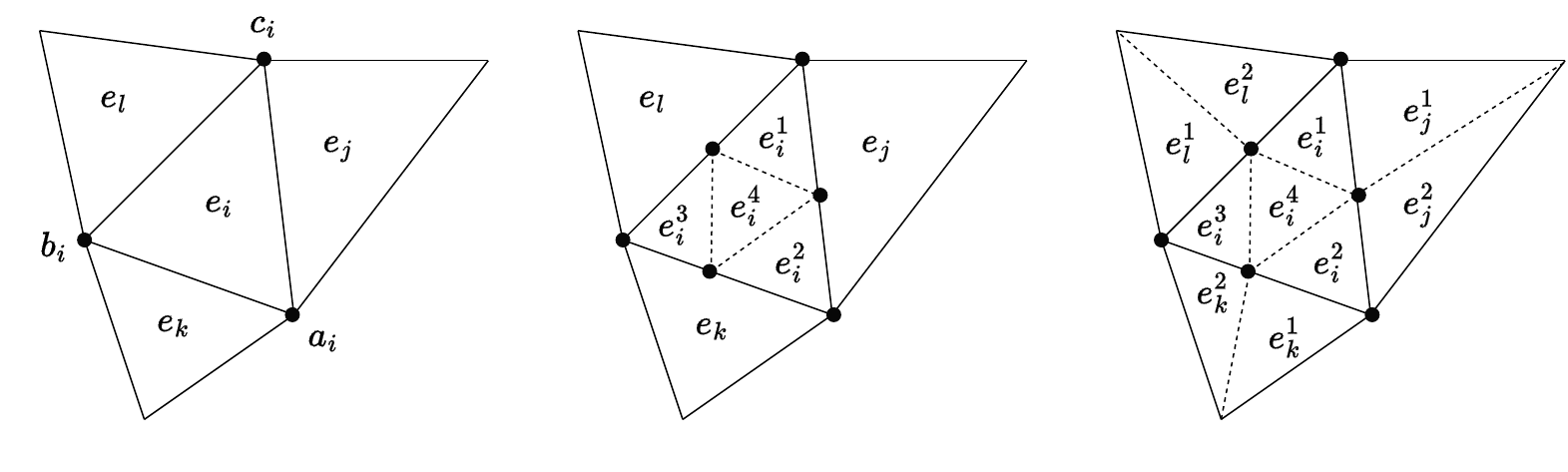}
	\caption{The red-green mesh refinement strategy \parencite{carstensen_adaptive_2004} for splitting element $e_i$. The element $e_i$, that satisfied the splitting criterion, and the three neighboring elements $e_i$, $e_j$, $e_k$ are replaced by 10 new elements. 
	}
	\label{fig:splitting}
\end{figure}

One benefit of the PINNs over the standard finite element interpolation is their ability to leverage all interpolation parameters to improve the accuracy of the interpolated solution. While the standard finite element method only solves for the values of the solution at the interpolation points, \emph{i.e.}, the mesh nodes, in the FENNI framework, the nodal coordinates can also be set as trainable. This increases the number of degrees of freedom to improve the interpolated solution.
The method
combines in a single step, solving the PDE and adapting the mesh with the target to find a solution minimizing the strain energy, for instance.
Conceptually, this allows getting closer to the continuous solution that actually minimizes the energy.
This optimization of nodal positions, known as \emph{r-adaptivity}, requires no additional steps beyond making the nodal coordinates trainable \AD{and has proven effective in \parencite{zhang2022hidenn}. We propose to extend the mesh adaptation to h-adaptivity based on the r-adaptivity dynamic. To that aim a new indicator is proposed to select the region that requires h-adaptivity (corresponding to the addition of nodes). The idea of this indicator is to rely on the movement of the nodes due to the r-adaptivity to drive the h-adaptivity.} 	\AD{Indeed,} r-adaptivity can highlight regions where mesh refinement is needed, as node movement toward certain areas requires higher density. To address this, an \emph{h-adaptivity} strategy driven by the movement of nodes from the r-adaptivity is proposed to improve the interpolation further. The areas where the mesh is locally refined by element splitting are defined as the regions where the elements shrink the fastest. We remark that the Jacobian determinant $J_i$, computed when evaluating the integral loss function, equals twice the area of element $e_i$. Therefore, we can conveniently define the element-splitting criterion directly with the Jacobian determinants. The smaller an element gets, the denser the zone should be. Thus, the mesh is refined in the area where nodes are moving, \emph{i.e.}, where the elements are getting smaller. The relative change of element area is determined as:
\begin{align}
	\Delta J_i^n = \frac{|J_i^{n-1}| - |J_i^n|}{|J_i^{n-1}|},
\end{align}
where $J_i^{n-1}$ and $J_i^n$ and the Jacobians determinant of element $e_i$ in $(n-1)$-th iteration and $n$-th iteration, respectively.
The elements satisfy the splitting criterion if $\Delta J_i^n > t_{\Delta J}$, where $t_{\Delta J}$ is a threshold value. An element $e_i$ that satisfies the splitting criterion is split using the red-green mesh refinement strategy \parencite{carstensen_adaptive_2004}. The refinement process, affecting the element $e_i$ and the neighbouring elements, is illustrated in Figure \ref{fig:splitting}. 

The maximum number of element splits is used as a parameter to determine the refinement level of the network. 
The elements that do not satisfy the splitting criterion but are affected by the splitting of a neighboring element are not considered to be split. Once an element reaches the maximal split value, it cannot be split but can be further deformed.

\subsection{Training} \label{sec:training_impl}

\subsubsection{Model construction}

First, a mesh of computational domain $\Omega$ with a given number of nodes is generated using Gmsh \parencite{gmsh_geuzaine2020three}. 
The model architecture is built based on the geometry and connectivity of the mesh. 
The Dirichlet boundary conditions are implemented by fixing the nodal values of boundary points \AD{to} the prescribed values. The nodal values are then made un-trainable.
The Neumann boundary conditions are incorporated into the loss function.

\subsubsection{Convergence criterion}
\AD{The model's parameters are optimized through an unsupervised training stage. In 1D, the training is performed by computing the loss on a set of uniformly sampled points in the domain $\Omega$. For the reference-element-based implementation, the sampled points' set directly } consists of the Gaussian quadrature points defined for each element. 
The training ends when either stagnation of loss or the maximal prescribed number of iterations is reached. The stagnation condition is met if the change of loss function between two consecutive iterations decreases below the defined threshold. The decrease in loss function is evaluated as:
\begin{align}
	\Delta l = 2 \frac{|l^n-l^{n-1}|}{|l^n + l^{n-1}|},
\end{align}
where $l^{n-1},\ l^n$ are the values of loss function in $n-1$- and $n$-th iteration.

\subsubsection{Multigrid training strategy} \label{sec:multi_level}

\AD{
	The model's interpretability enables a multigrid training strategy to be designed, thus improving the training efficiency. All trainable parameters can indeed be directly reused to initialize subsequent models. In the multigrid training approach, a model is initially trained on a coarse mesh, and its parameters are then used to initialize a model on a finer mesh. This process can be iterated across increasingly refined meshes until the target resolution is reached.}

Let us consider a domain $\Omega$ discretized by the meshes with $np_1$ and $np_2$ nodes, where $np_2>np_1$. 
Model $M_{np_1}$ is trained until reaching convergence. Afterward, the nodal values associated with the mesh nodes $X_{np_2}$ of the new model $M_{np_2}$ 

is initialized by evaluating the initial model $M_{np_1}(X_{np_2})$.
With the reference-element-based implementation, the element IDS needs to be provided when evaluating the model $M_{np_1}$. If a node of the finer mesh is not located in any element of the coarse mesh, we evaluate the model $M_{np_1}$ in the closest node within the original mesh. This generally happens for the boundary nodes.

\section{Numerical results} \label{sec:res}
In this section, we investigate the accuracy of the proposed method. On 1D and 2D linear elasticity problems, we demonstrate the effect of the choice of integration method used to evaluate the integral loss function, the loss function, and the mesh adaptivity. 

The implementation of the method \parencite{daby-seesaramNeuROM2024b} is freely available online\footnote{\url{https://github.com/AlexandreDabySeesaram/NeuROM}}.
The results presented in this section can be reproduced using interactive demos \parencite{skardovaFENNIIPaperDemo2024}. The Jupyter notebooks for the demos can be run online\footnote{\url{https://alexandredabyseesaram.github.io/FENNI-I-paper-demo/}}.

\subsection{Considered problems}

\KS{\subsubsection{1D: Bar }}
We use a simple 1D example, a modification of a problem used in \parencite{zhang2021hierarchical} where one of the zero Dirichlet boundary conditions is replaced by a non-zero boundary value. We consider
a beam with both ends fixed, under a body force $b(x)$:
\begin{align}
	&\frac{\mathrm{d}}{\mathrm{d}x}\left(AE \frac{\mathrm{d}u}{\mathrm{d}x}\right) + b\left(x\right) = 0 \quad \mathrm{in}\ \Omega = [0,L], \label{eq:1D_problem_1}\\
	&u\left(0\right) = 0,\ u\left(L\right) = u_L, \label{eq:1D_problem_2}
\end{align}
where $u(x)$ is the displacement field, $A=\SI{1}{mm\squared}$ is the beam section area, $E = 175$ MPa is the beam stiffness and $u_L = \SI{5e-4}{mm}$ is the right Dirichlet boundary condition. The body force $b(x)$ is prescribed as:
\begin{align}
	b\left(x\right) = -\frac{\left(4 \pi^2 \left(x - x_1\right)^2 - 2 \pi\right)}{e^{\left(\pi (x - x_1)^2\right)}} - \frac{\left(8 \pi^2 \left(x - x_2\right)^2 - 4 \pi\right)}{e^{\left(\pi \left(x - x_2\right)^2\right)}},
\end{align} 
where 
\begin{equation}
	\begin{cases}
		x_1 = \SI{2.5}{mm} \\
		x_2 = \SI{7.5}{mm} 
	\end{cases}
\end{equation}
The analytical solution for the displacement and its derivative is the following:
\begin{align}
	\nonumber
	u_{a}\left(x\right) = &\frac{1}{AE}\left(e^{-\pi\left((x-x_1\right)^2} - e^{-x_1^2\pi}\right) +
	\frac{2}{AE}\left(e^{-\pi\left(x-x_2\right)^2} - e^{-x_2^2\pi}\right)\\
	&- \frac{1}{10AE}\left(e^{-x_1^2\pi} - e^{-x_2^2\pi}\right)x + u_L x,\\
	\nonumber
	\frac{\mathrm{d}u_{a}}{\mathrm{d}x} = &\frac{1}{AE}\left(\pi e^{-\pi(x-x_1)^2(x-x_1)}\right) +
	\frac{4}{AE}\left(-\pi e^{-\pi(x-x_2)^2} (x-x_2)\right)\\
	&- \frac{1}{10AE}\left(e^{-x_1^2\pi} - e^{-x_2^2\pi}\right) + u_L
\end{align}

\KS{\subsubsection{2D: Plate with a hole}}

\begin{figure}[H]
	\centering
	\includegraphics[width=0.6\linewidth]{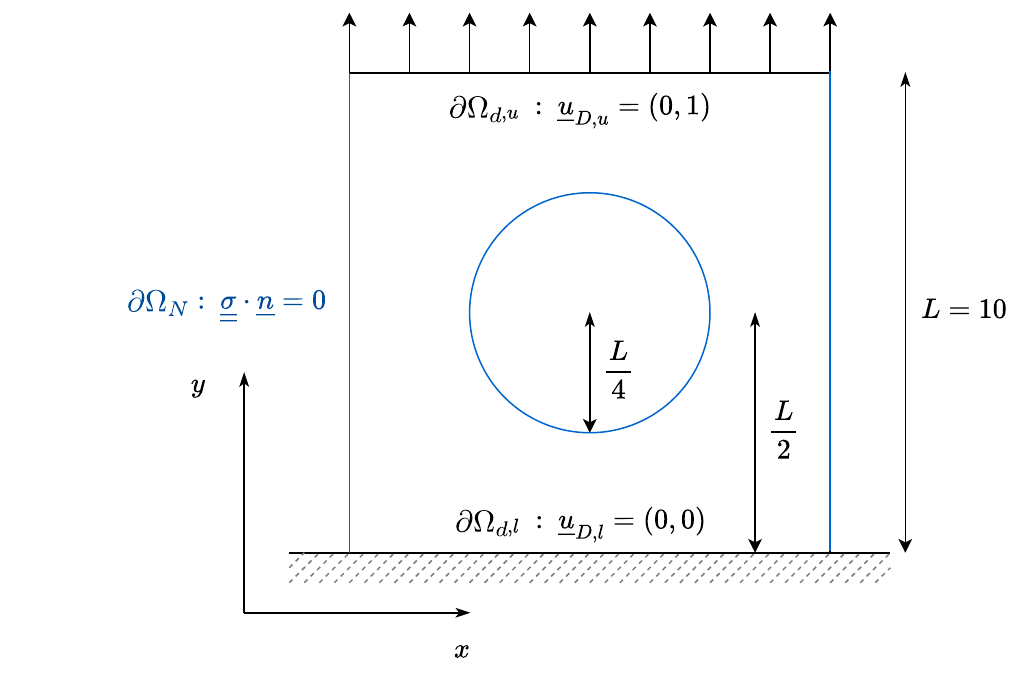} 
	\caption{2D problem setting.}
	\label{fig:2D_problem}
\end{figure}

The effect of r- and rh-adaptivity are further illustrated in 2D. Considering small deformations, the governing equations can be written as:
\begin{align}
	-\nabla \cdot \doubleunderline{\sigma} &= 0 \quad \mathrm{in}\ \Omega,\\
	\doubleunderline{\sigma} &= \lambda \mathrm{tr}(\doubleunderline{\varepsilon})\mathbb{I} + 2\mu \doubleunderline{\varepsilon} \quad \mathrm{in}\ \Omega,\\
	\doubleunderline{\varepsilon} &= \frac{1}{2}\left(\nabla \underline{u} + \left( \nabla \underline{u} \right)^T \right) \quad \mathrm{in}\ \Omega,\\
	\doubleunderline{\sigma} \cdot \underline n &= 0 \quad \mathrm{on}\  \AD{\dSn}
	,\\
	\underline u &= \underline u_D \quad \mathrm{on}\ \AD{\dSd},
\end{align}
where $\lambda = 1.25, \mu = 1.0$ are the Lamé coefficients. The domain $\Omega$ is a bar with one hole shown in Figure \ref{fig:2D_problem}.

The displacement $\underline u$, obtained as the output of the model, is used to compute the stress tensor $\doubleunderline{\sigma}$ and von Mises stress $\sigma^r_{VM}$, defined as $\sigma_{VM} = \sqrt{\frac{3}{2}\doubleunderline{\sigma}^{\mathrm{dev}}:\doubleunderline{\sigma}^{\mathrm{dev}}}$, where $\doubleunderline{\sigma}^{\mathrm{dev}}$ is the stress deviator tensor $\doubleunderline{\sigma}^{\mathrm{dev}} = \doubleunderline{\sigma} - \frac{1}{3}\mathrm{tr}(\doubleunderline{\sigma})\mathbb{I}$. 

\KS{\subsubsection{Error metrics}}

The error of the model results with respect to the analytical solution in 1D is evaluated using the following norm:
\begin{align}
	\|e_u\|_2 = \frac{\left( \sum_{i=1}^{N_e} \left(u\left(x_i\right) - u_a\left(x_i\right)\right)^2\right)^{1/2}}{\left(\sum_{i=1}^{N_e} u_a\left(x_i\right)^2\right)^{1/2}},
	\label{eq:error_u}
\end{align}
where the evaluation points $x_i$, are sampled uniformly within the interval $\Omega$. The number of evaluation points is $N_e = 1000$.
A similar norm is used to evaluate the error in strain:
\begin{align}
	\|e_{\nabla u}\|_2 = \frac{\left( \sum_{i=1}^{N_e} \left(\frac{\mathrm{d}u}{\mathrm{d}x}\left(x_i\right) - \frac{\mathrm{d}u_{a}}{\mathrm{d}x}\left(x_i\right)\right)^2\right)^{1/2}}{\left(\sum_{i=1}^{N_e} \frac{\mathrm{d}u_{a}}{\mathrm{d}x}\left(x_i\right)^2\right)^{1/2}},
	\label{eq:error_du}
\end{align}

In 2D, a numerical solution computed on a very fine mesh (37 988 elements, 19 368 nodes) using the finite element method is used as a reference solution. 
We compare the displacement and von Mises obtained by the proposed method with the reference values. The normalized error norms $\|e_{u_x}\|_2$, $\|e_{u_y}\|_2$, $\|e_{\sigma_{VM}}\|_2$ are defined analogically to the norms used in 1D (\ref{eq:error_u}). For example:
\begin{align}
	\|e_{u_x}\|_2 = \frac{\left( \sum_{i} (u^r_x\left(x_i\right) - u_x\left(x_i\right))^2\right)^{1/2}}{\left(\sum_{i=1}^{N_e} u^r_x\left(x_i\right)^2\right)^{1/2}},
\end{align}
where the evaluation points $x_i$ are defined as the centres of the elements of the reference mesh, \emph{i.e.}, $N_e$ = 37 988.
Besides the global norms, we also use one local metric -- the value of maximal von Mises stress. The relative maximal von Mises stress is computed as a ratio between the predicted value and the reference value:
\begin{align}
	\sigma^{max}_{VM} = \frac{\max_{i=1:N_e} \sigma_{VM}\left(x_i\right)}{\max_{i=1:N_e} \sigma^r_{VM}\left(x_i\right)}
\end{align}

\KS{\subsection{Impact of the choice of the integration method}
	In this section, we aim to evaluate the impact of the choice of the integration method -- trapezoidal rule and Gaussian quadrature rule -- on the accuracy of the results. This comparison is conducted using a 1D example and potential energy loss function for fixed and r-adaptive mesh\AD{es}. The nodal values were initialized uniformly by value $0.1$. The training was performed using the PyTorch implementation of the L-BFGS algorithm with a strong Wolfe line search.}

The minimized loss, representing the system's potential energy (\ref{eq:loss_general_potential}) is defined for the governing equations (\ref{eq:1D_problem_1}), (\ref{eq:1D_problem_2}) as follows:
\begin{align}
	J(u) = \frac{1}{2} \int_\Omega AE \left(\frac{\mathrm{d}u}{\mathrm{d}x}\right)^2 \mathrm{d}x - \int_\Omega ub \ \mathrm{d}x
	\label{eq:loss_potential_en_1D}
\end{align}

To study the convergence to the analytical solution, the domain $\Omega$ is uniformly discretized by 10, 21, 41, 80, 160, and 324 nodes, respectively. 
With the trapezoidal rule, the training was performed using three datasets containing 10, 20, and 30 sampling points per element. The sampling points are distributed uniformly. 
With the Gaussian quadrature rule, the training was performed using 2, 3, 4, and 5 integration points. 

\begin{figure}[h]
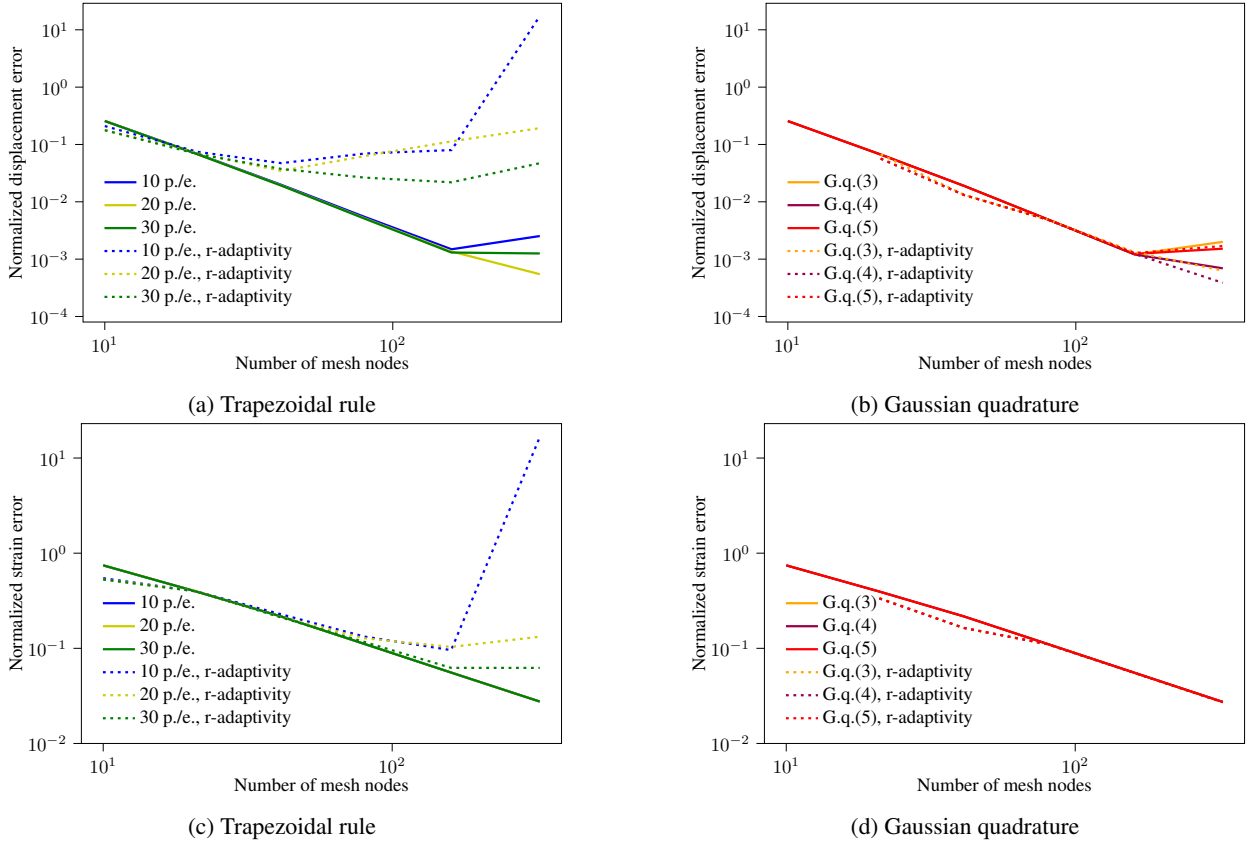

	\centering
	\begin{subfigure}{0.45\textwidth}
		\includegraphics[width=\textwidth]{1D_adapt_trapezoid_u.tikz}
		\caption{Trapezoidal rule}
	\end{subfigure}
	\hfill
	\begin{subfigure}{0.45\textwidth}
		\includegraphics[width=\textwidth]{1D_adapt_Gauss_u.tikz}
		\caption{Gaussian quadrature}
	\end{subfigure}\\
	\begin{subfigure}{0.45\textwidth}
		\includegraphics[width=\textwidth]{1D_adapt_trapezoid_grad.tikz}
		\caption{Trapezoidal rule}
	\end{subfigure}
	\hfill
	\begin{subfigure}{0.45\textwidth}
		\includegraphics[width=\textwidth]{1D_adapt_Gauss_grad.tikz}
		\caption{Gaussian quadrature}
	\end{subfigure}
	\caption{Normalized displacement (a,b) and strain (c,d) error. Two methods for evaluating the potential energy loss functions were used: the trapezoidal rule (a,c) and the Gaussian quadrature rule (b,d). For the trapezoidal rule, the notation $n$ p./e. indicate the number of samples per element. The sampling points are distributed uniformly and are not moved during the training. For the Gaussian quadrature rule, the notation G.q.($n$) indicates the use of the $n$-point Gaussian quadrature rule.  }
	\label{fig:1D_adapt_convergence}
\end{figure}

The normalized displacement and strain error for fixed and r-adaptive mesh are shown in Figure \ref{fig:1D_adapt_convergence}. Note that the strain is obtained by differentiating the model output $u$ with respect to the model input $x$. 
For the fixed-mesh training, the accuracy when using the trapezoidal rule and the Gaussian quadrature rule are comparable. 
Compared to fixed mesh, the accuracy is improved only for very coarse mesh (10 - 21 nodes). For finer meshes, the error with adaptive mesh is significantly higher than with the fixed mesh. 
For the coarsest mesh (10 nodes), the training with r-adaptive mesh diverges for all tested numbers of quadrature points. For the other mesh resolutions, the accuracy of displacement and strain is the same or higher compared to the fixed mesh. 

\begin{figure}[h]
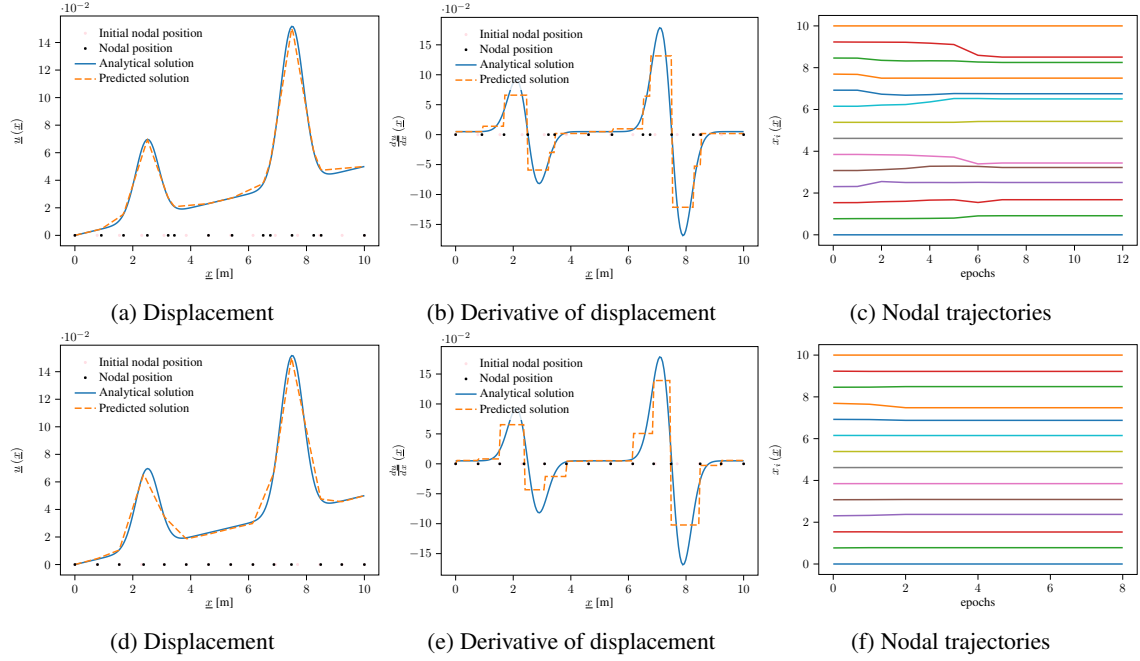

	\centering
	\begin{subfigure}{0.3\textwidth}
		\includegraphics[width=\textwidth]{Gauss_Displacement.tikz}
		\caption{Displacement}
	\end{subfigure}
	\begin{subfigure}{0.3\textwidth}
		\includegraphics[width=\textwidth]{Gauss_Gradient.tikz}
		\caption{Derivative of displacement}
	\end{subfigure}
	\begin{subfigure}{0.3\textwidth}
		\includegraphics[width=\textwidth]{Gauss_Trajectories.tikz}
		\caption{Nodal trajectories}
	\end{subfigure}
	
	\begin{subfigure}{0.3\textwidth}
		\includegraphics[width=\textwidth]{Trapezoid_Displacement.tikz}
		\caption{Displacement}
	\end{subfigure}
	\begin{subfigure}{0.3\textwidth}
		\includegraphics[width=\textwidth]{Trapezoid_Gradient.tikz}
		\caption{Derivative of displacement}
	\end{subfigure}
	\begin{subfigure}{0.3\textwidth}
		\includegraphics[width=\textwidth]{Trapezoid_Trajectories.tikz}
		\caption{Nodal trajectories}
	\end{subfigure}
	\caption{The results on 14-nodes mesh obtained using the 5-point Gaussian quadrature (a-c) and trapezoidal rule with $30$ sampling points per element (d-f). }
	\label{fig:1D_solution}
\end{figure}

Figure \ref{fig:1D_solution} illustrates the difference between the two integration methods on r-adaptive mesh with 14 nodes. The figure shows that the use of the Gaussian quadrature rule leads to a more accurate displacement and strain due to more efficient adjustment of the nodal positions.

\AD{It is worth noting}, that when using the trapezoidal rule with r-adaptive mesh, the training points are not redistributed during the training. On the other hand, with the Gaussian quadrature rule, quadrature points are assigned to the deformed elements as the mesh is adapted. \AD{This aspect, combined with the absence of integration error allowed by the Gauss quadrature, probably explains its use as being more robust than integration by the trapezoid method. In the remainder of this article, the Gauss quadrature will be employed to compute integral terms in the loss functions.}

\subsection{Impact of the choice of the loss function}
\KS{In this section, the results obtained with the potential energy loss function are compared with the residual and weak formulation loss function. This comparison is performed using the 1D example and fixed mesh. The nodal values were initialized uniformly to $0.1$, and the L-BFGS optimizer was used. }

\subsubsection{Weak formulation loss function}

The weak formulation loss function defined by equation (\ref{eq:loss_general_weak}) has the following form:
\begin{align}
	J(u) = \sum_{i=1}^{np} \left( \int_{\Omega_i} AE \frac{\mathrm{d}u}{\mathrm{d}x} \frac{\mathrm{d}u^*_i}{\mathrm{d}x} \mathrm{d}x - \int_{\Omega_i} u^*_ib \mathrm{d}x \right)^2 \ \mathrm{d}x,
	\label{eq:loss_weak_eq_en_1D}
\end{align}
where $u^*_i$ is the $i$-th \AD{test} function and ${\Omega_i}$ is its support.

The comparison with the potential energy loss function is shown in Figure \ref{fig:1D_convergence_weak_eq}.
For the meshes with 10 - 80 nodes, the results are identical to the results obtained using the potential energy loss. 

As the number of mesh nodes increases, the training slows down prohibitively, which argues in favor of using the potential energy-based loss in the remainder of this article.

\subsubsection{Residual loss}

For problems the governing equations (\ref{eq:1D_problem_1}),(\ref{eq:1D_problem_2}), the minimized loss consists of the following PDE residual term and the compatibility term:
\begin{align}
	J(u) = \lambda_1\mathrm{MSE}\left( \frac{\mathrm{d} {u}_x}{\mathrm{d}x} - \frac{AE}{b}\right) + \lambda_2\mathrm{MSE}\left( \frac{\mathrm{d}u}{\mathrm{d}x} - {u}_x\right)
	\label{eq:loss_mixed_1D}
\end{align}
where $\frac{\mathrm{d}u}{\mathrm{d}x}$ is the strain computed based on predicted displacement $u$, an ${u}_x$ is the predicted strain. The weights are set to $\lambda_1 = L$, $\lambda_2 = 1$.
To ensure the same order of $\frac{\mathrm{d}u}{\mathrm{d}x}$ and ${u}_x$ in the compatibility term, quadratic shape functions are used in the model for prediction of $u$ and linear shape functions for prediction of ${u}_x$. 
Training was conducted using uniformly distributed \AD{sample points}, with 10, 25, 50, or 75 points per element.

\begin{figure}[h]
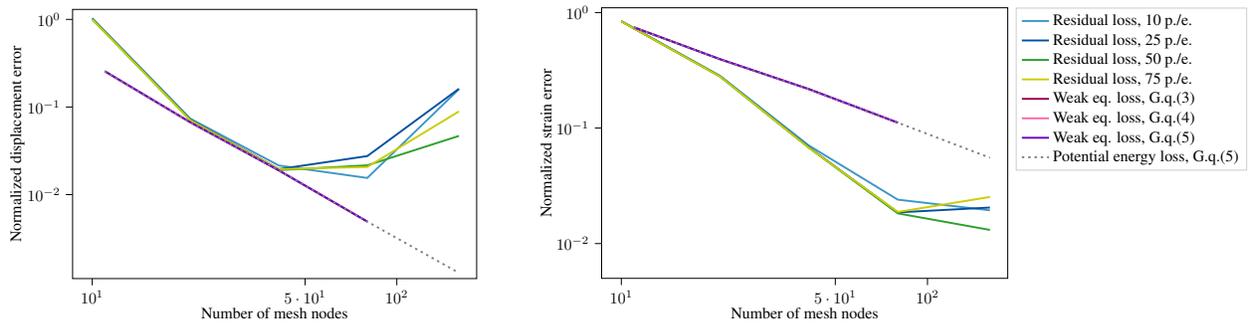

	\centering
	\begin{subfigure}{0.38\textwidth}
		\includegraphics[width=\textwidth]{1D_mixed_form_u}
	\end{subfigure}
	\hfill
	\begin{subfigure}{0.573\textwidth}
		\includegraphics[width=\textwidth]{1D_mixed_form_grad.tikz}
	\end{subfigure}
	\caption{Comparison of residual loss, weak formulation loss, and potential energy loss. For the finest tested mesh (160 nodes), the training with weak formulation loss failed to converge in the prescribed maximal number of iterations.}
	\label{fig:1D_convergence_weak_eq}
\end{figure}

The comparison with the potential energy loss function is shown in Figure \ref{fig:1D_convergence_weak_eq}.
The normalized strain error is lower for all mesh resolutions. The normalized displacement error is comparable for medium-coarse meshes but significantly higher for coarse and fine meshes.

We remark that the model with potential energy loss only predicts displacement $u$, represented by linear shape functions. The strain, computed by differentiating the displacement, is, therefore, only piece-wise constant. The difference between the two methods shown in Figure \ref{fig:1D_convergence_weak_eq} is partly due to this difference in the interpolation order.

\KS{\subsection{Impact of the choice of optimizers and training strategies}
	In this section, we investigate the effect of the choice of optimizer -- Adam with learning rate = $10^{-3}$ and L-BFGS. The comparison is done on the 2D example problem using the potential-energy loss function, which was shown to be the most efficient, as shown above. The nodal values were initialized uniformly by value $0.5$. }

\begin{figure}
	\begin{minipage}[b]{.46\linewidth}
		\centering
		\includegraphics[width=\linewidth]{2D_Adam_LBFGS.tikz}
		\captionof{figure}{Convergence to reference solution for training performed using Adam and L-BFGS.}
		\label{fig:2D_convergence}
	\end{minipage}
	\hfill
	\begin{minipage}[b]{.46\linewidth}
		\centering
		\begin{tabular}{ *{3}{c} }
			44 nodes &	144	nodes & 484 nodes \\ \hline
			\includegraphics[width=0.3\linewidth, trim={27cm 2cm 27cm 0cm},clip]{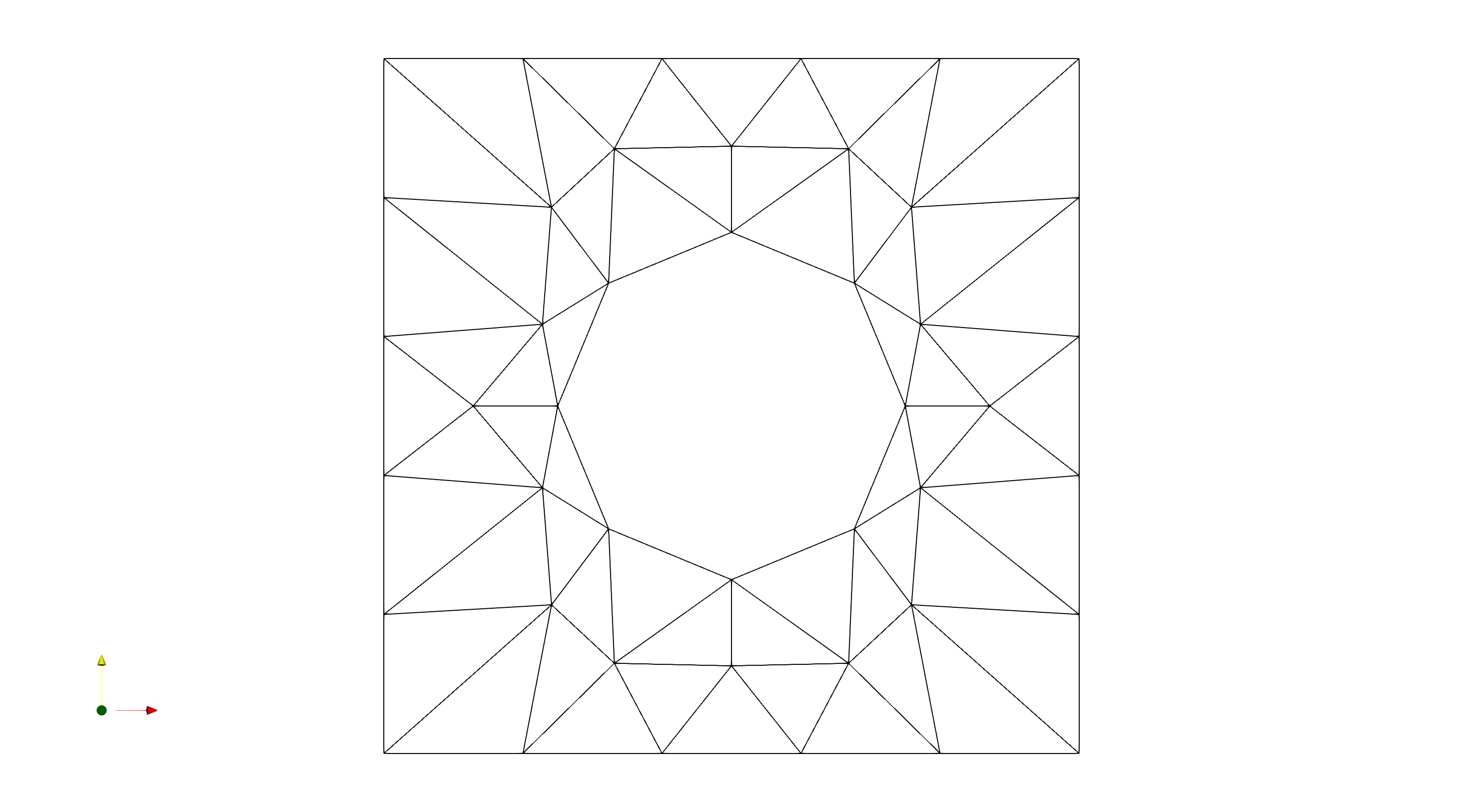} &
			\includegraphics[width=0.3\linewidth, trim={27cm 2cm 27cm 0cm},clip]{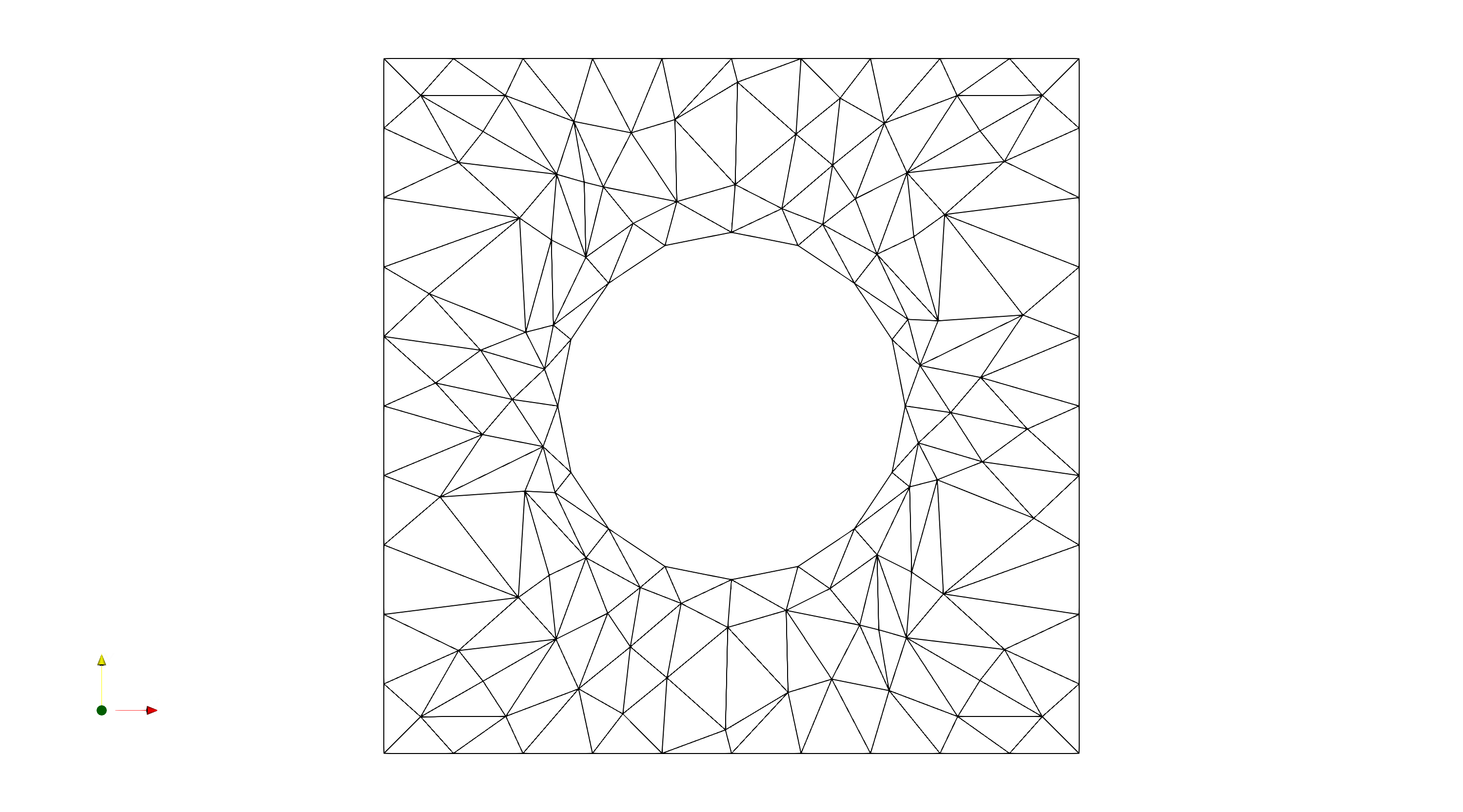} &
			\includegraphics[width=0.3\linewidth, trim={27cm 2cm 27cm 0cm},clip]{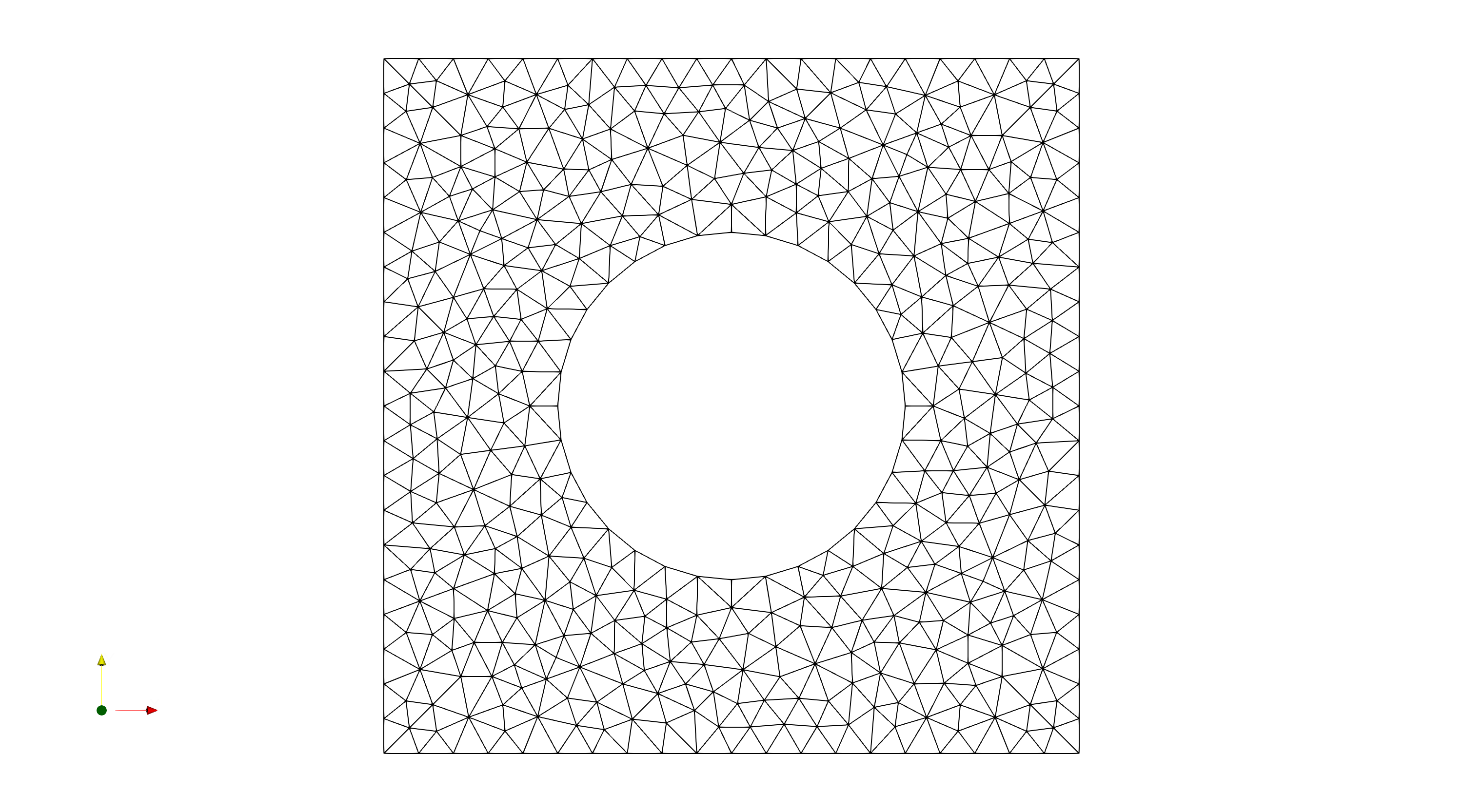} \\
		\end{tabular}
		\begin{tabular}{ *{2}{c} }
			1804 nodes &	7088 nodes \\ \hline
			\includegraphics[width=0.3\linewidth, trim={27cm 2cm 27cm 0cm},clip]{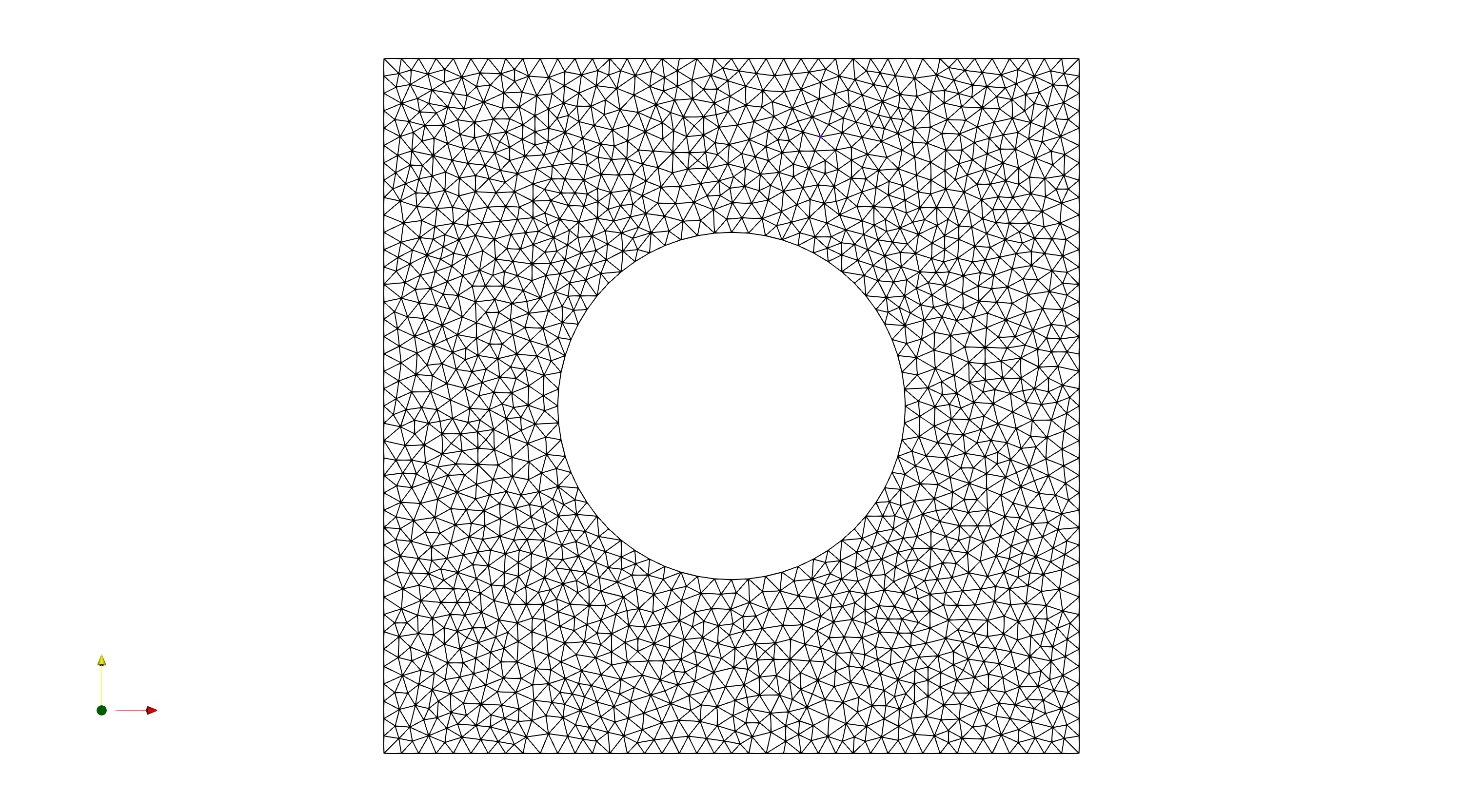} &
			\includegraphics[width=0.3\linewidth, trim={27cm 2cm 27cm 0cm},clip]{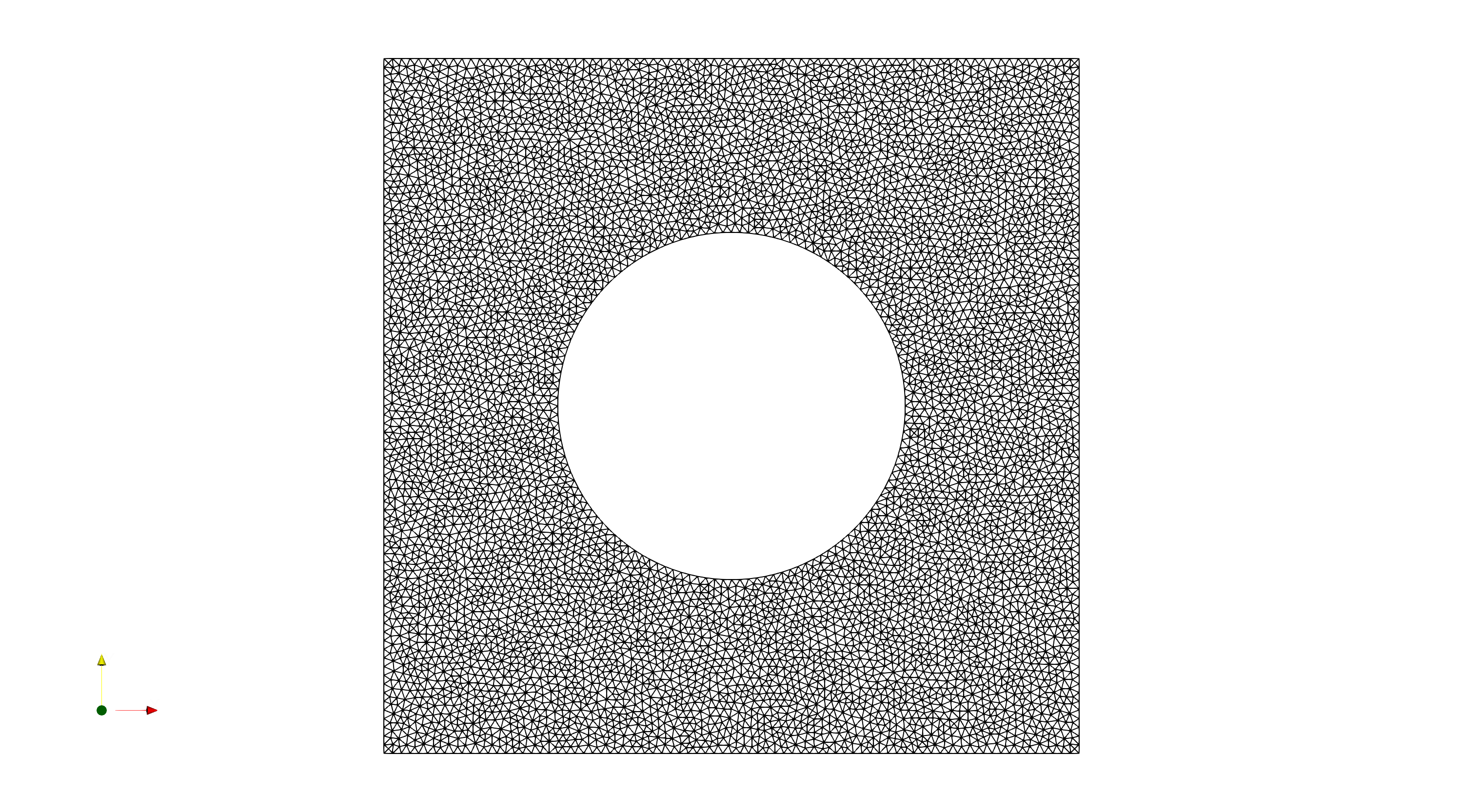} 
		\end{tabular}
		\captionof{table}{Discretization of domain $\Omega$.}
		\label{tab:2D_meshes}
	\end{minipage}
\end{figure}

The minimized loss:
\begin{align}
	J(u) = \frac{1}{2} \int_\Omega \doubleunderline{\sigma}:\doubleunderline{\varepsilon} \ \mathrm{d}x,
	\label{eq:loss_potential_en_2D}
\end{align}
is evaluated using the Gauss quadrature rule with 1 integration point. 
To study the convergence to the reference solution, the domain $\Omega$ is discretized by triangular meshes with 44,	144, 484, 1 804, and 7 088 nodes. The meshes are shown in Table \ref{tab:2D_meshes}. 
The normalized error in displacement and stress for computations on fixed mesh are shown in Figure \ref{fig:2D_convergence}. For the coarser meshes, the accuracy of both optimizers is comparable while L-BFGS provides significantly better results for finer meshes.

\begin{figure}[h]
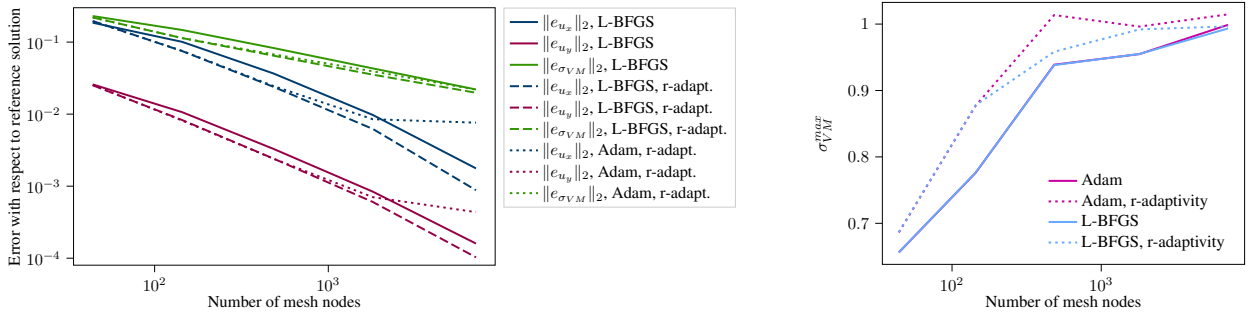

	\centering
	\begin{subfigure}{0.59\textwidth}
		\includegraphics[width=\textwidth]{2D_adapt_multi_Adam_LBFGS.tikz}
		\caption{The effect of r-adaptivity on the norm of error in displacement components and the von Mises stress. }
	\end{subfigure}
	\hfill
	\begin{subfigure}{0.35\textwidth}
		\includegraphics[width=\textwidth]{2D_max_stress.tikz}
		\caption{The maximal value of von Mises stress obtained with and without r-adaptivity. }
	\end{subfigure}
	\caption{The effect of r-adaptivity with increasing mesh resolution. A multigrid training strategy with scaling factor 2 was used.}
	\label{fig:2D_r}    
\end{figure}

\begin{figure}[h]
	\centering
	\begin{subfigure}{0.3\textwidth}
		\centering
		\includegraphics[height=\textwidth, trim={28cm 2cm 22cm 2cm},clip]{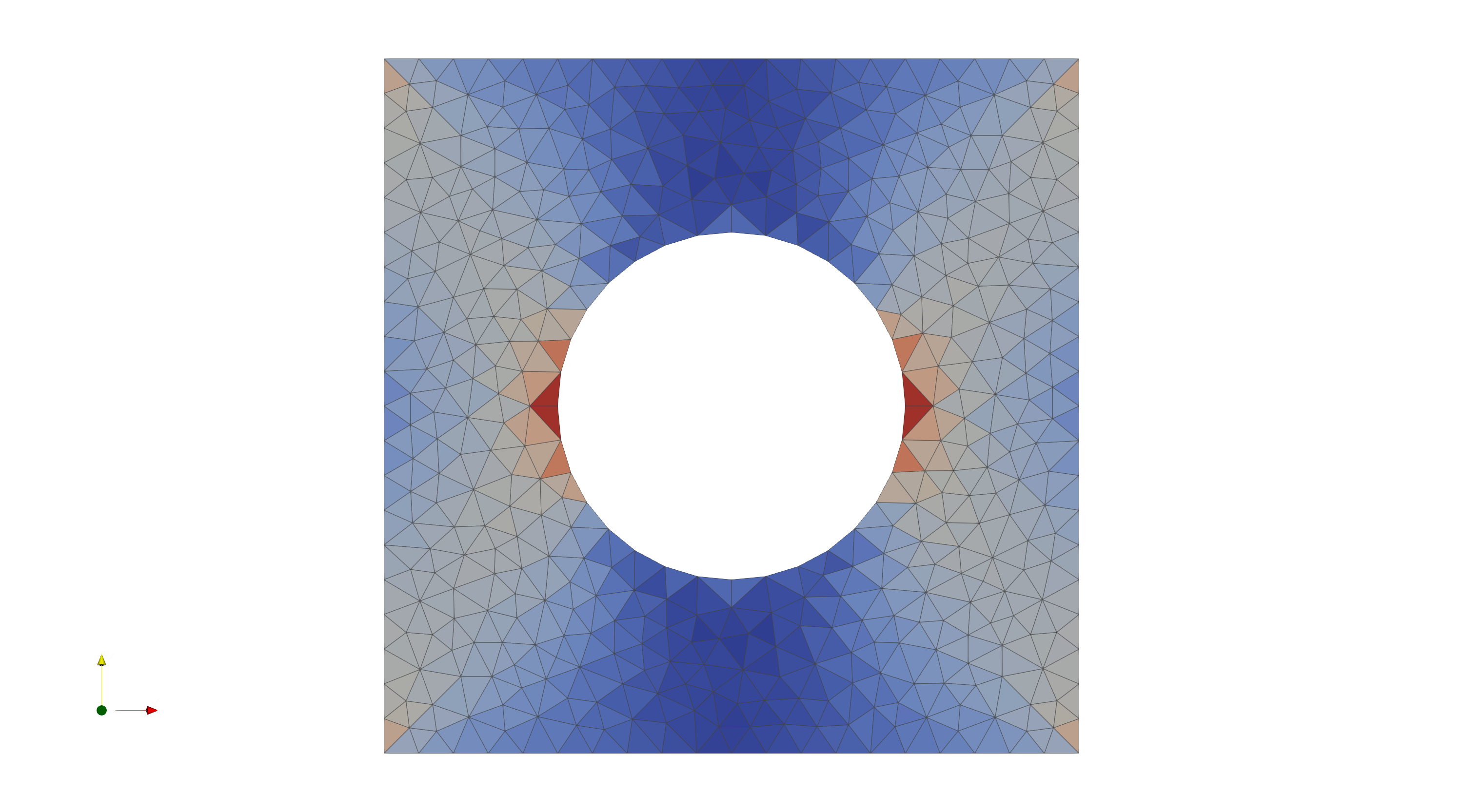}
		\caption{Fixed mesh}
	\end{subfigure}
	\begin{subfigure}{0.3\textwidth}
		\centering
		\includegraphics[height=\textwidth, trim={28cm 2cm 22cm 2cm},clip]{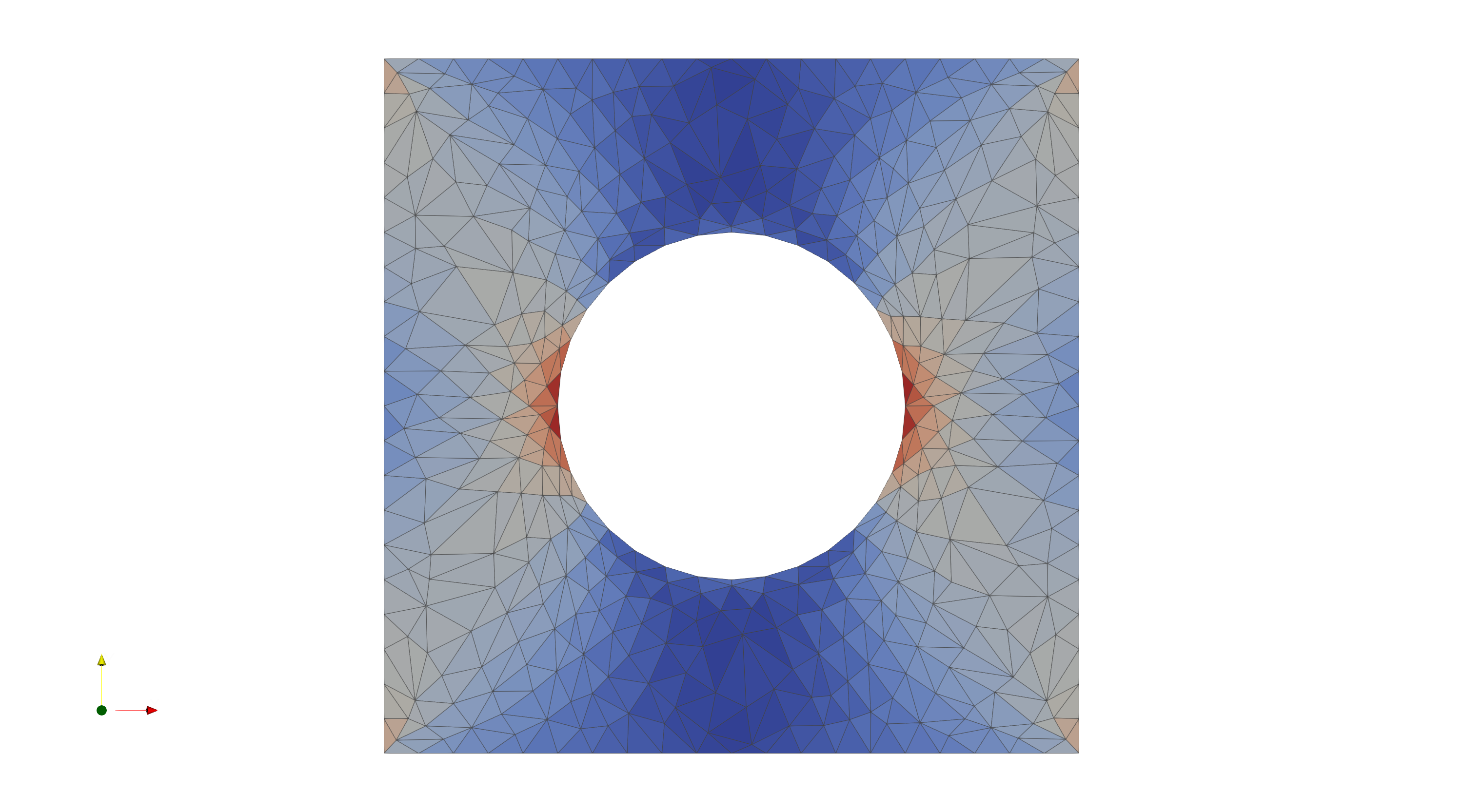}
		\caption{r-adaptivity}
	\end{subfigure}
	\begin{subfigure}{0.3\textwidth}
		\centering
		\includegraphics[height=\textwidth, trim={28cm 2cm 11cm 2cm},clip]{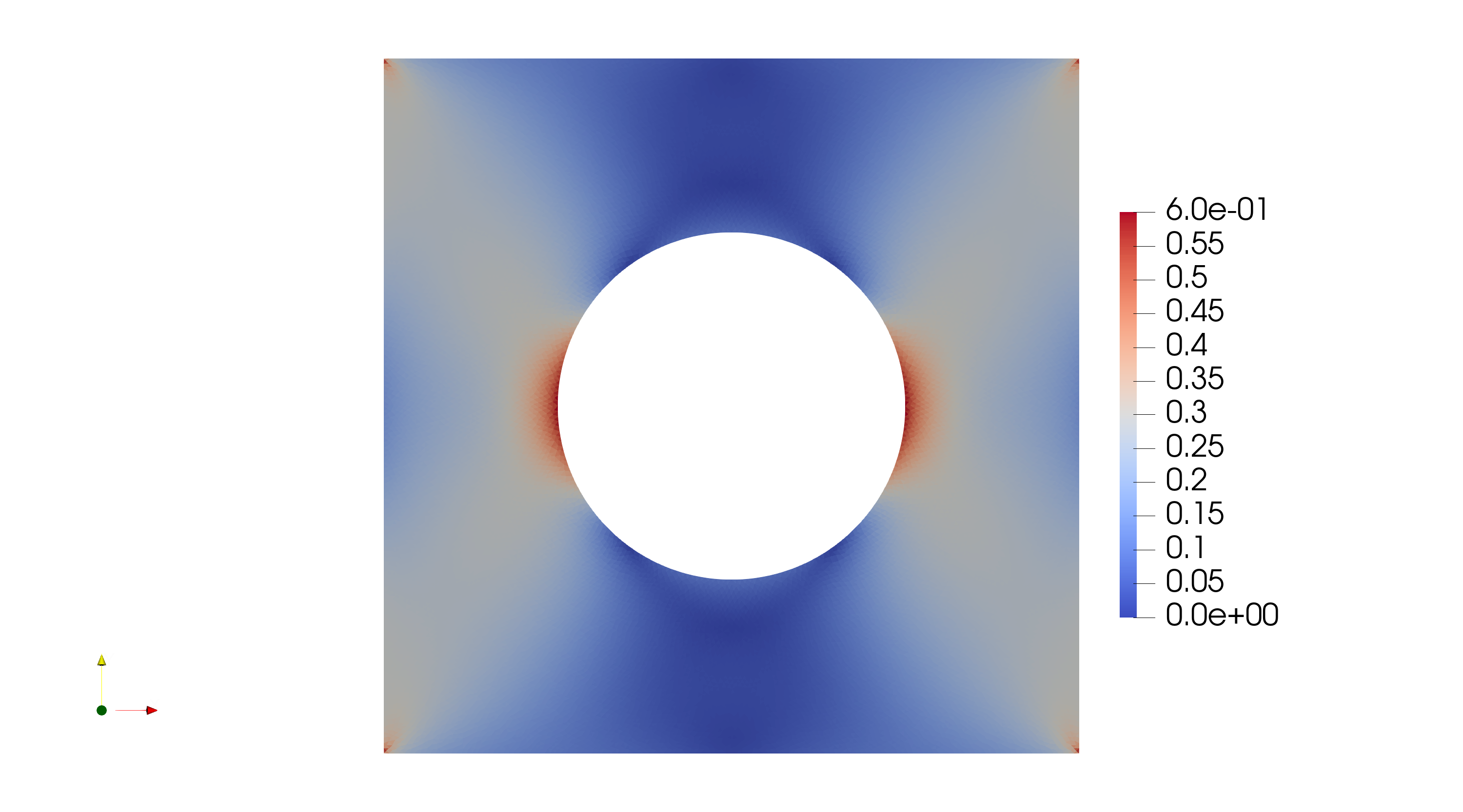}
		\caption{Reference solution}
	\end{subfigure}
	\caption{The reference von Mises stress $\sigma^r_{VM}$ and $\sigma_{VM}$ computed based on displacement obtained from the NN model with fixed mesh and r-adaptivity. The L-BFGS optimizer was used in both cases. }
	\label{fig:2D_r_solution}    
\end{figure}

The normalized displacement and stress errors and the relative maximal von Misses stress for the r-adaptive meshes are shown in Figure \ref{fig:2D_r}, where they are compared with results obtained on the fixed mesh. The use of L-BFGS optimizer with r-adaptive mesh leads to higher accuracy in all metrics. The Adam optimizer was shown to be the same or less accurate in all predictions, except the maximal von Misses stress for the fine mesh resultuions. 
The von Mises stress on the fixed and r-adaptive mesh is visualized in Figure \ref{fig:2D_r_solution}. The figure shows that the mesh gets refined in the areas with high-stress values.

The training with r-adaptive mesh was performed using the multigrid training strategy, described in Section \ref{sec:multi_level}. The training starts on the coarsest mesh, and after reaching convergence, the mesh is refined, and the training continues. The sequence of mesh resolutions is defined by scaling factor 2.

The \AD{computational speed-up enabled by the multigrid training strategy} is provided in \cref{tab:2D_r_adapt_bis}. 
With Adam optimizer, the single-level training didn't converge in the prescribed maximal number of iterations for the finest mesh (7088 nodes). With L-BFGS, the single-level training failed for two finest meshes (1804 and 7088 nodes). The experiment shows that for finer meshes, the multigrid approach is not only faster but, in some cases, necessary to reach convergence.

\begin{table}[h!]
	\centering
	\begin{tabular}{ c c c c c c} 
		& Scale 1 & Scale 2 & Scale 3 & Scale 4 & Scale 5 \\
		& 44 nodes & 	144	 nodes & 484 nodes &	1 804 nodes & 7 088 nodes \\
		\hline
		Adam  & 1 & 0.48	& 0.8 &	3.1 & $\infty$\\ 
		L-BFGS & 1 & 1 &  0.78 & $\infty$ & $\infty$
	\end{tabular}
	\caption{\AD{Computational Speed-up of the multigrid training with scaling factor 2, with respect to the single-level training. }}
	\label{tab:2D_r_adapt_bis}
\end{table}

\begin{figure}
	\begin{minipage}[t][][t]{.47\linewidth}
		\centering
		\begin{minipage}[t]{\textwidth}
			\centering
			\includegraphics[width=\textwidth]{2D_rh_adapt_multi.tikz} 
			\vspace{0.3cm}
		\end{minipage}
		\begin{minipage}[t]{\textwidth}
			\centering
			\includegraphics[width=\textwidth]{2D_rh_adapt_max_stress.tikz}
			\caption{The effect of rh-adaptivity on the accuracy of displacement and von Mises stress. }
			\label{fig:2D_rh_norms}
			\vspace{2cm}
		\end{minipage}
		
	\end{minipage}\hfill
	\begin{minipage}[t]{.45\linewidth}
		\centering
		\begin{tabular}{ *{3}{c} }
			Initial mesh &	r-adaptivity & rh-adaptivity \\
			44 nodes & 44 nodes & 74 nodes \\\hline
			\includegraphics[width=0.3\linewidth, trim={27cm 2cm 27cm 0cm},clip]{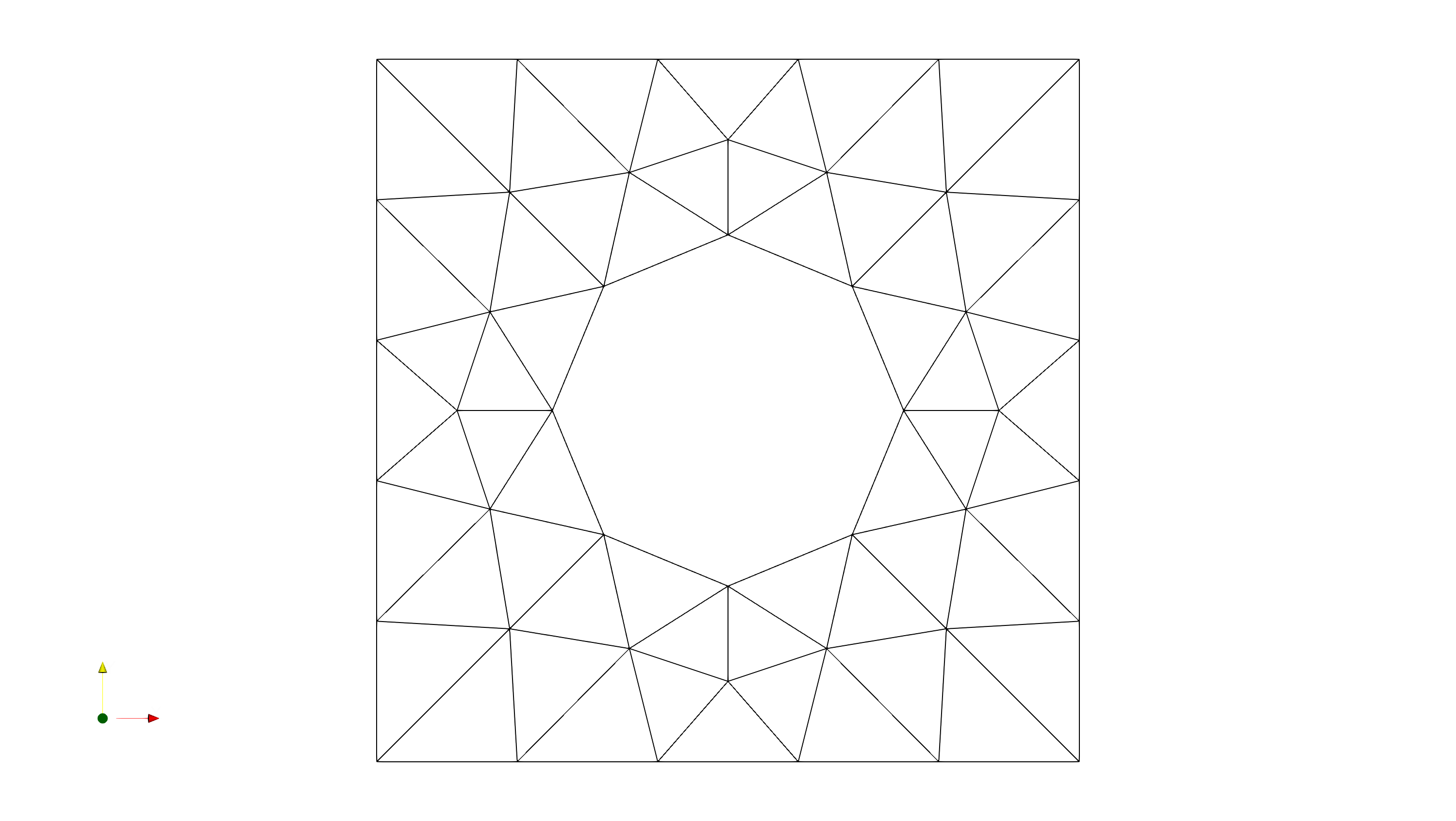} &
			\includegraphics[width=0.3\linewidth, trim={27cm 2cm 27cm 0cm},clip]{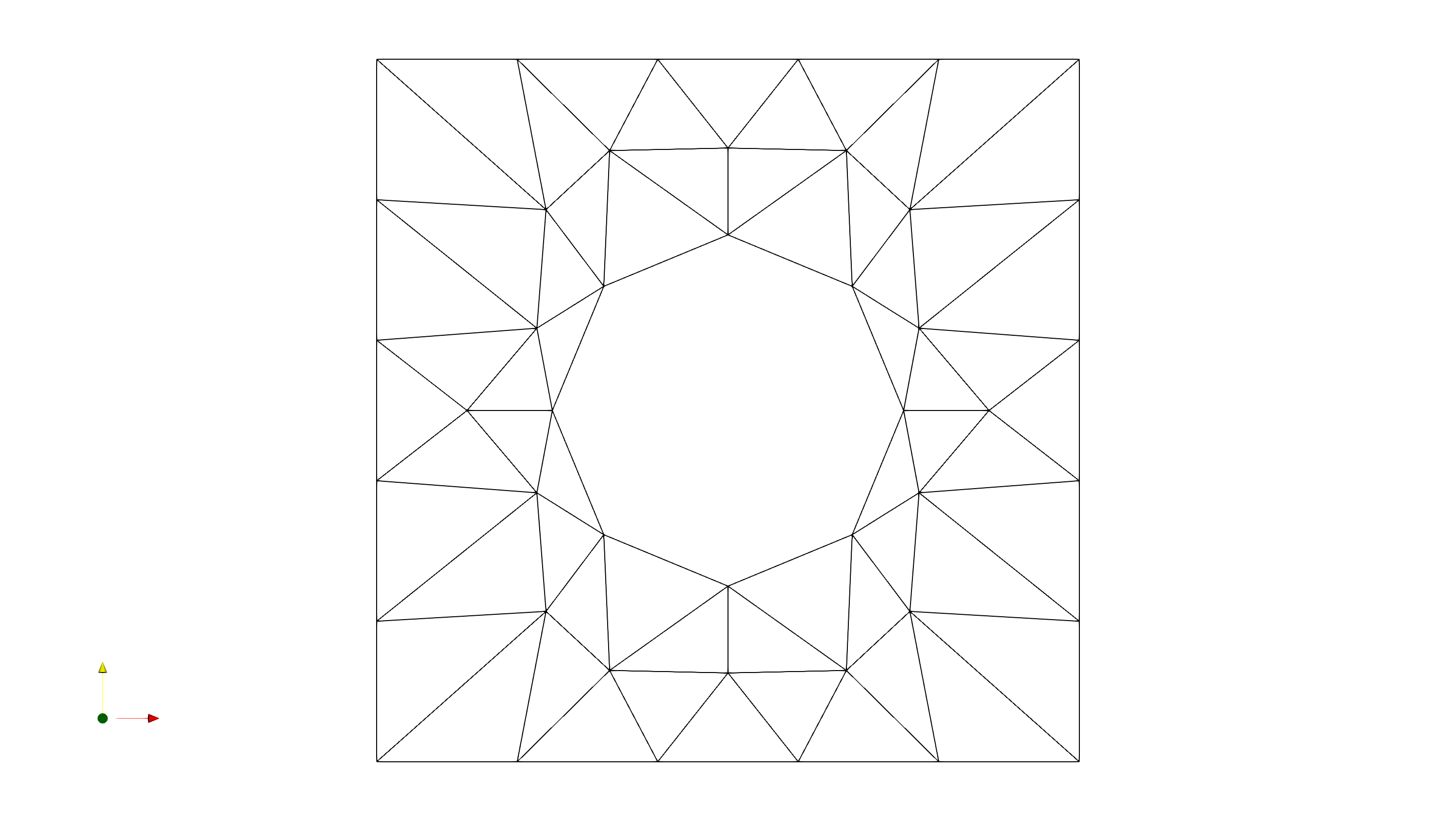} &
			\includegraphics[width=0.3\linewidth, trim={27cm 2cm 27cm 0cm},clip]{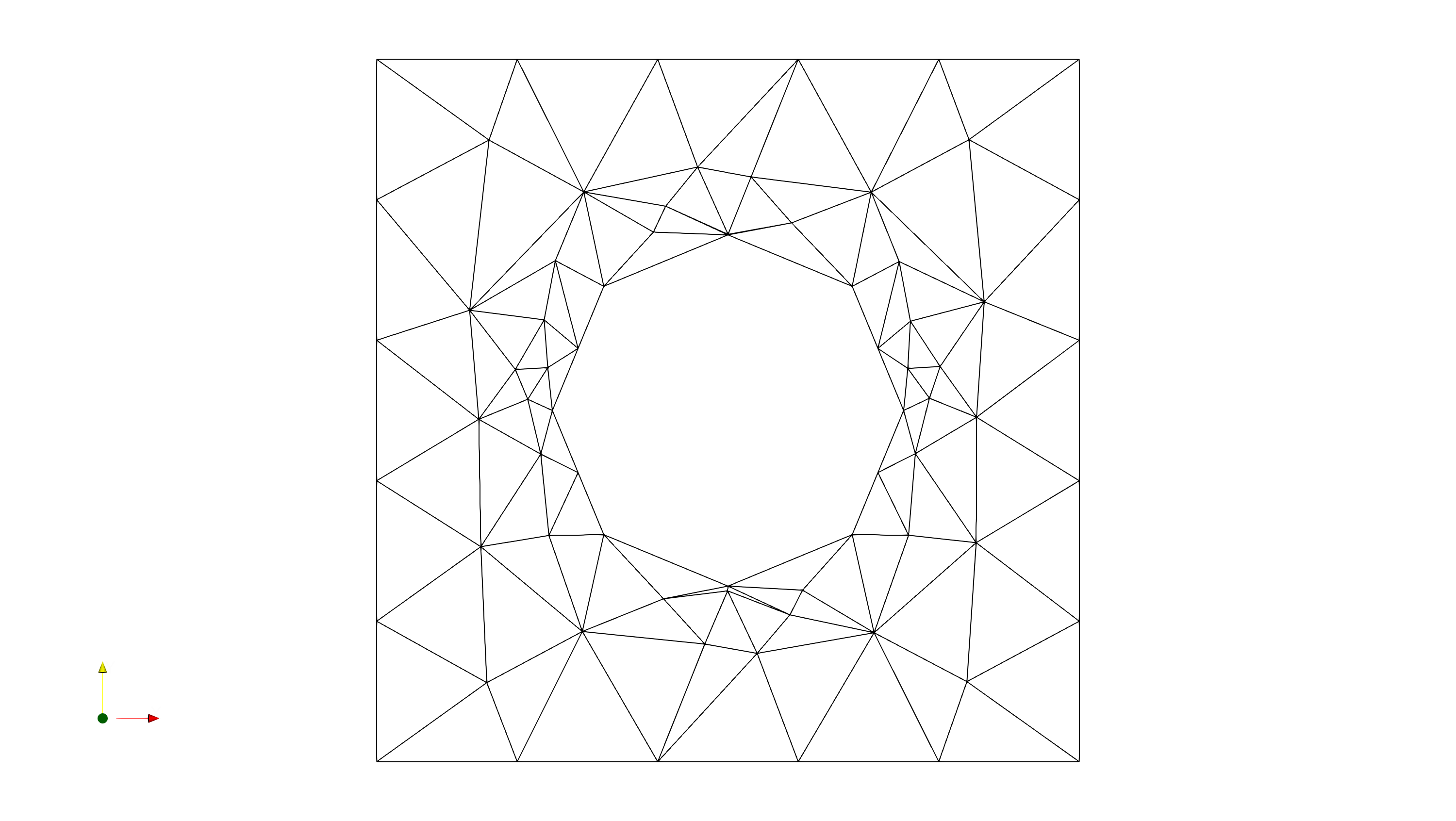} \\
			144 nodes & 144 nodes & 292 nodes \\\hline
			\includegraphics[width=0.3\linewidth, trim={27cm 2cm 27cm 0cm},clip]{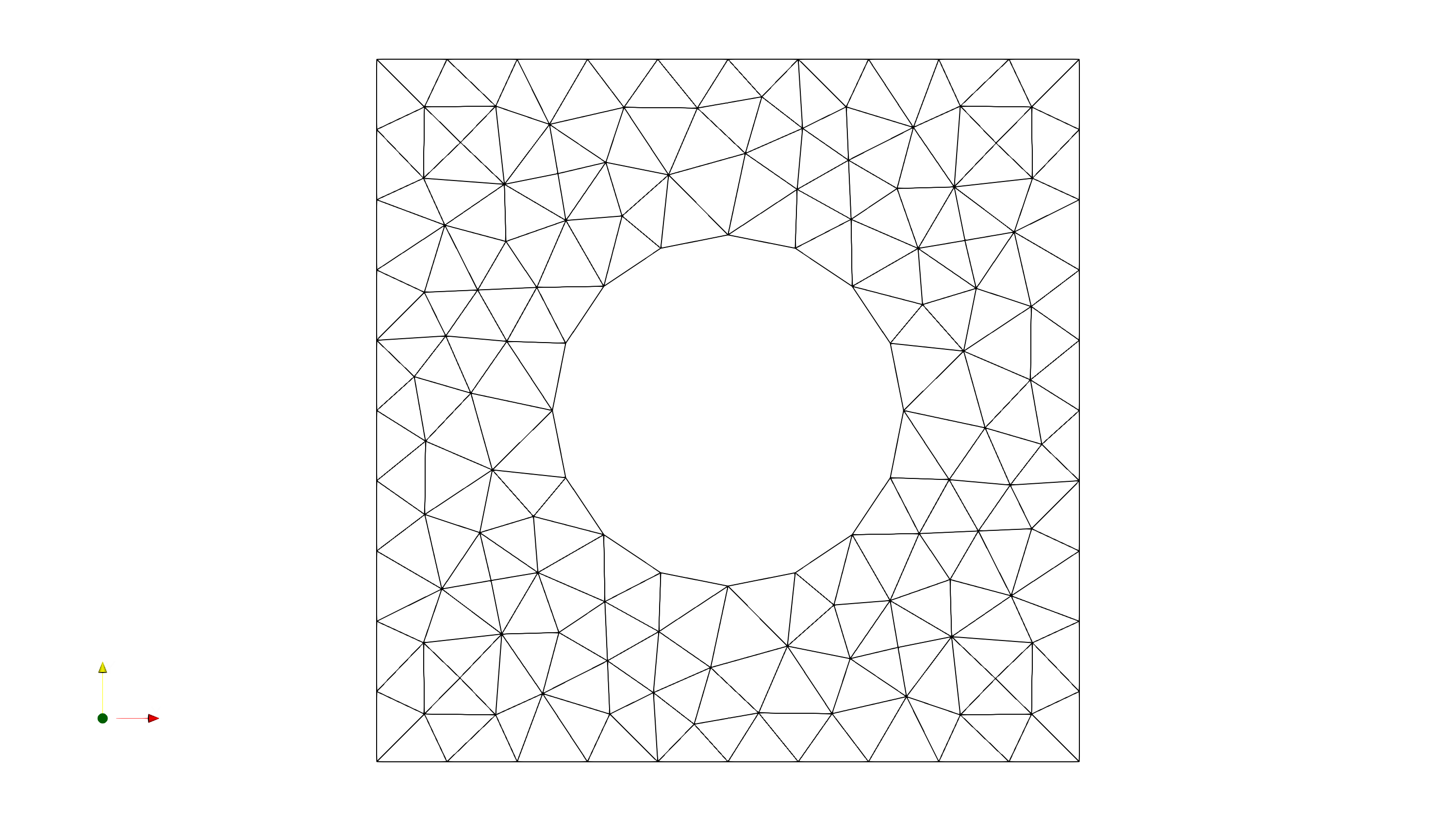} &
			\includegraphics[width=0.3\linewidth, trim={27cm 2cm 27cm 0cm},clip]{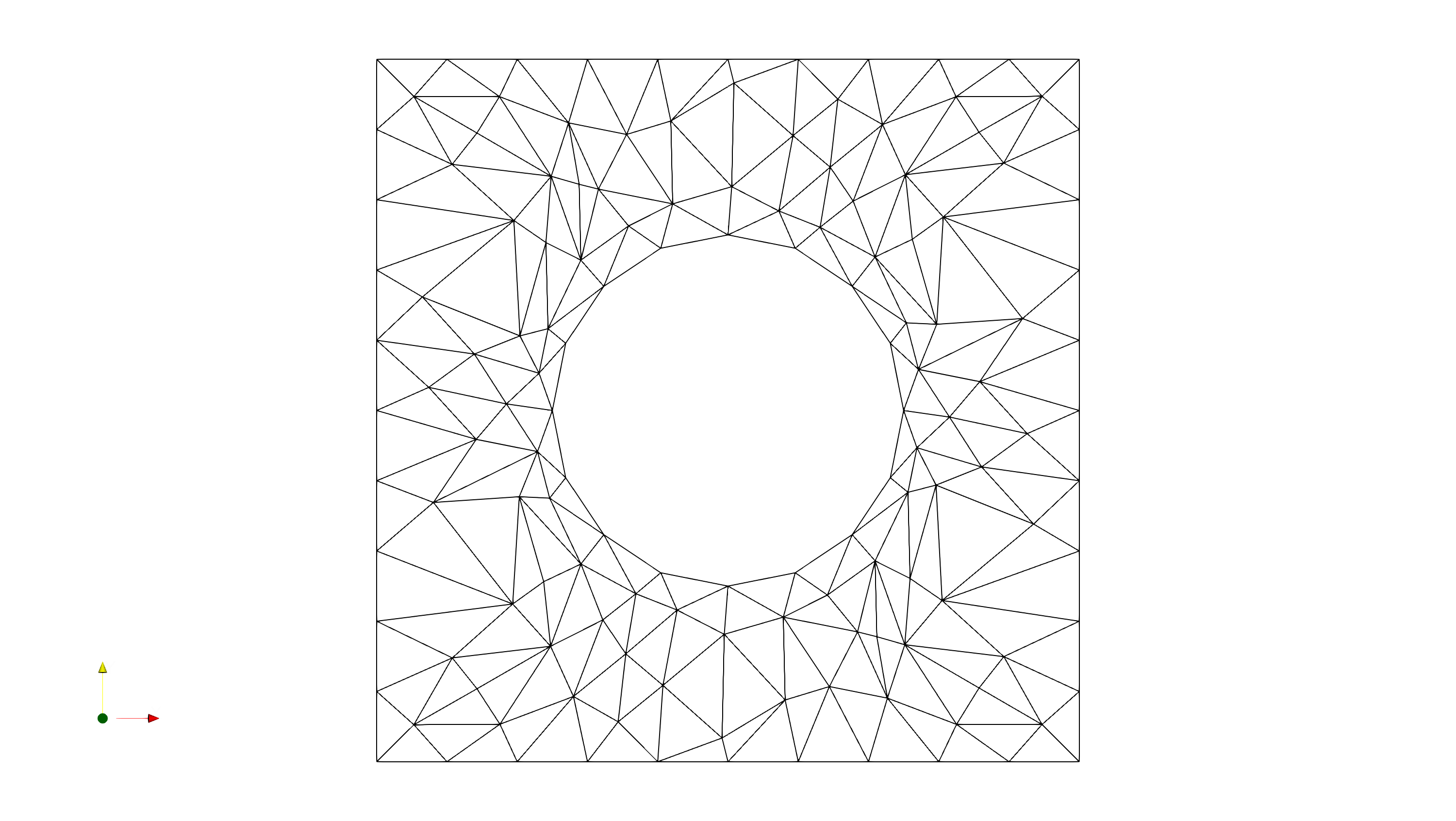} &
			\includegraphics[width=0.3\linewidth, trim={27cm 2cm 27cm 0cm},clip]{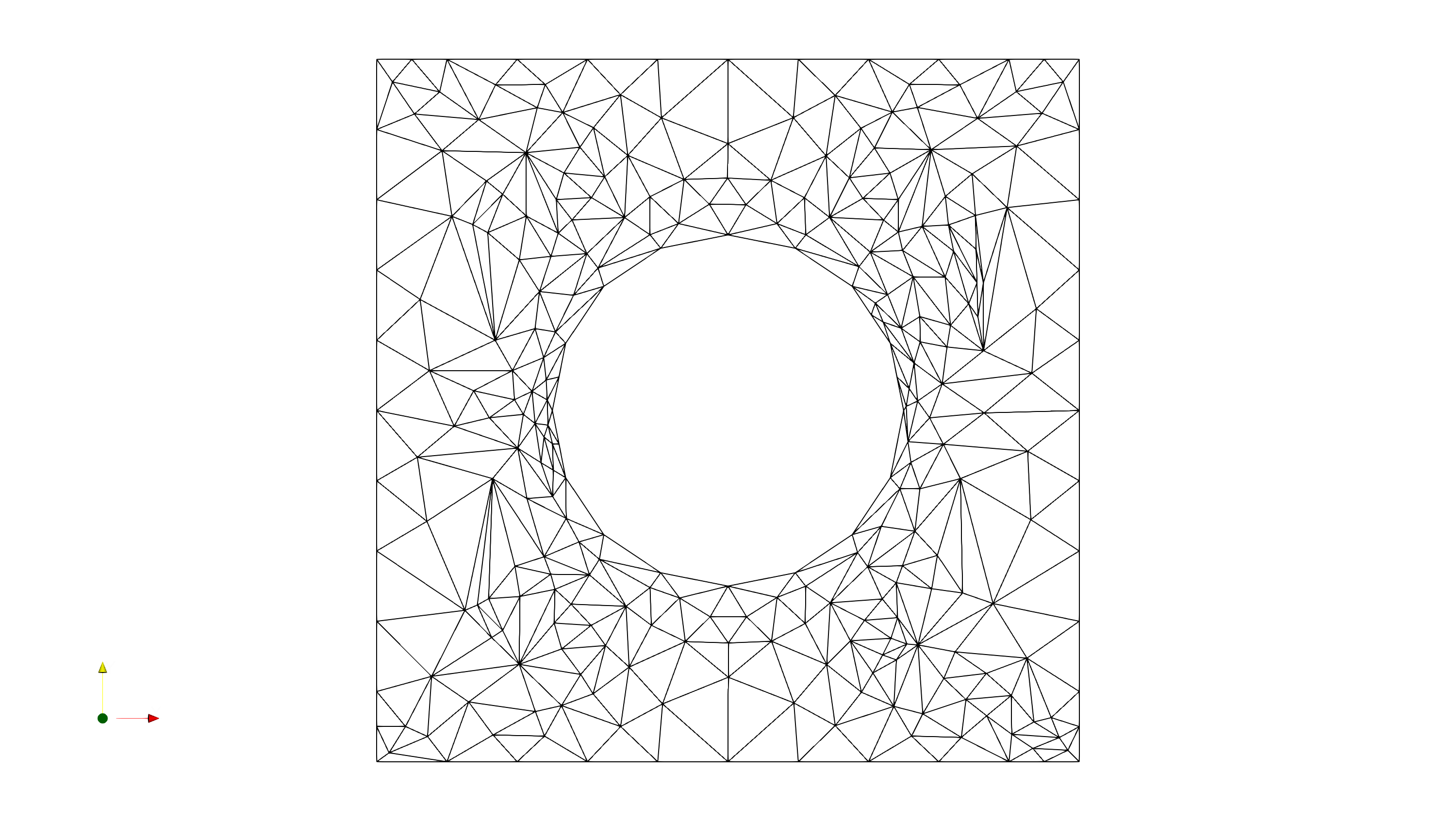} \\
			484 nodes & 484 nodes & 684 nodes \\\hline
			\includegraphics[width=0.3\linewidth, trim={27cm 2cm 27cm 0cm},clip]{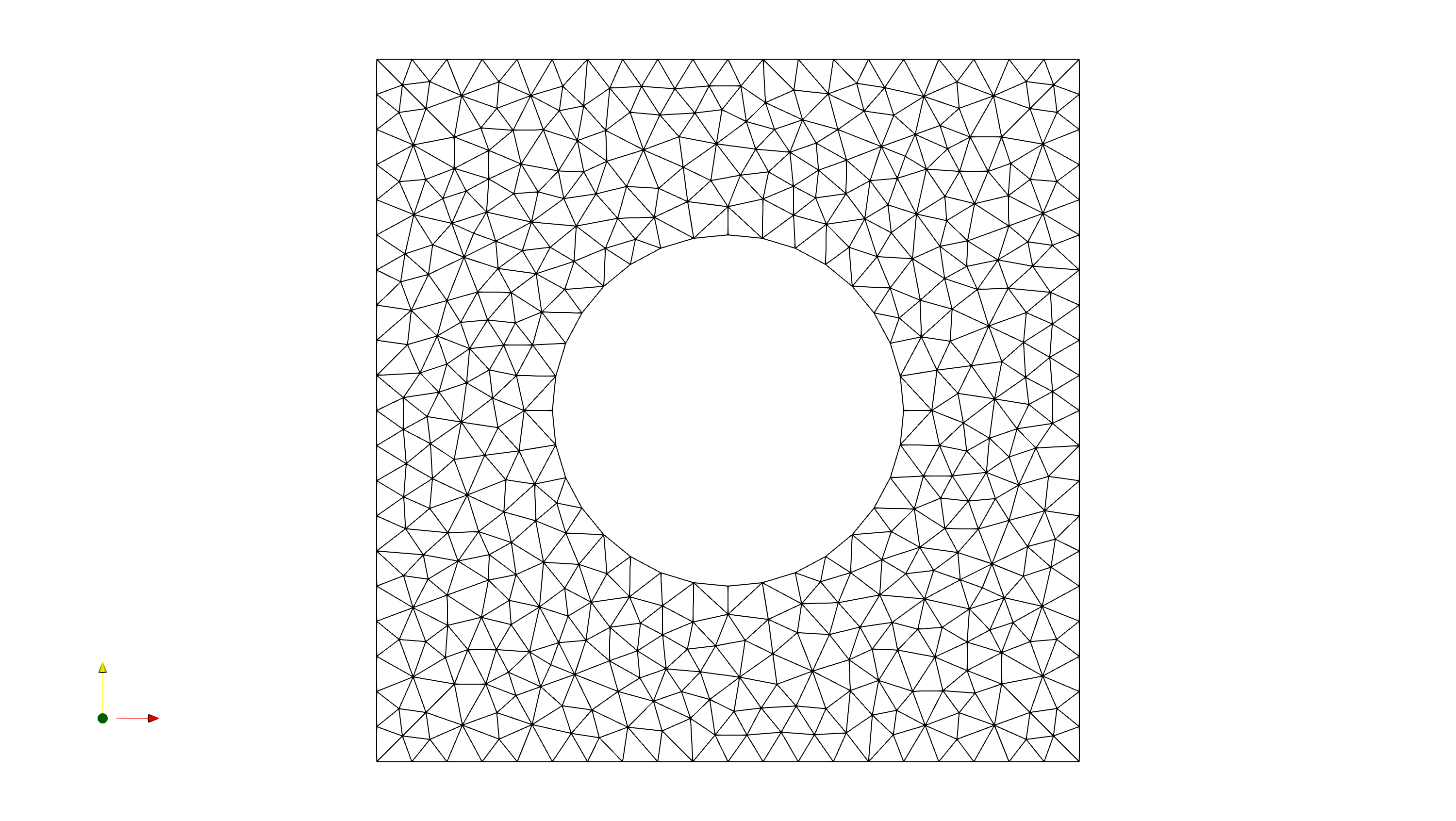} &
			\includegraphics[width=0.3\linewidth, trim={27cm 2cm 27cm 0cm},clip]{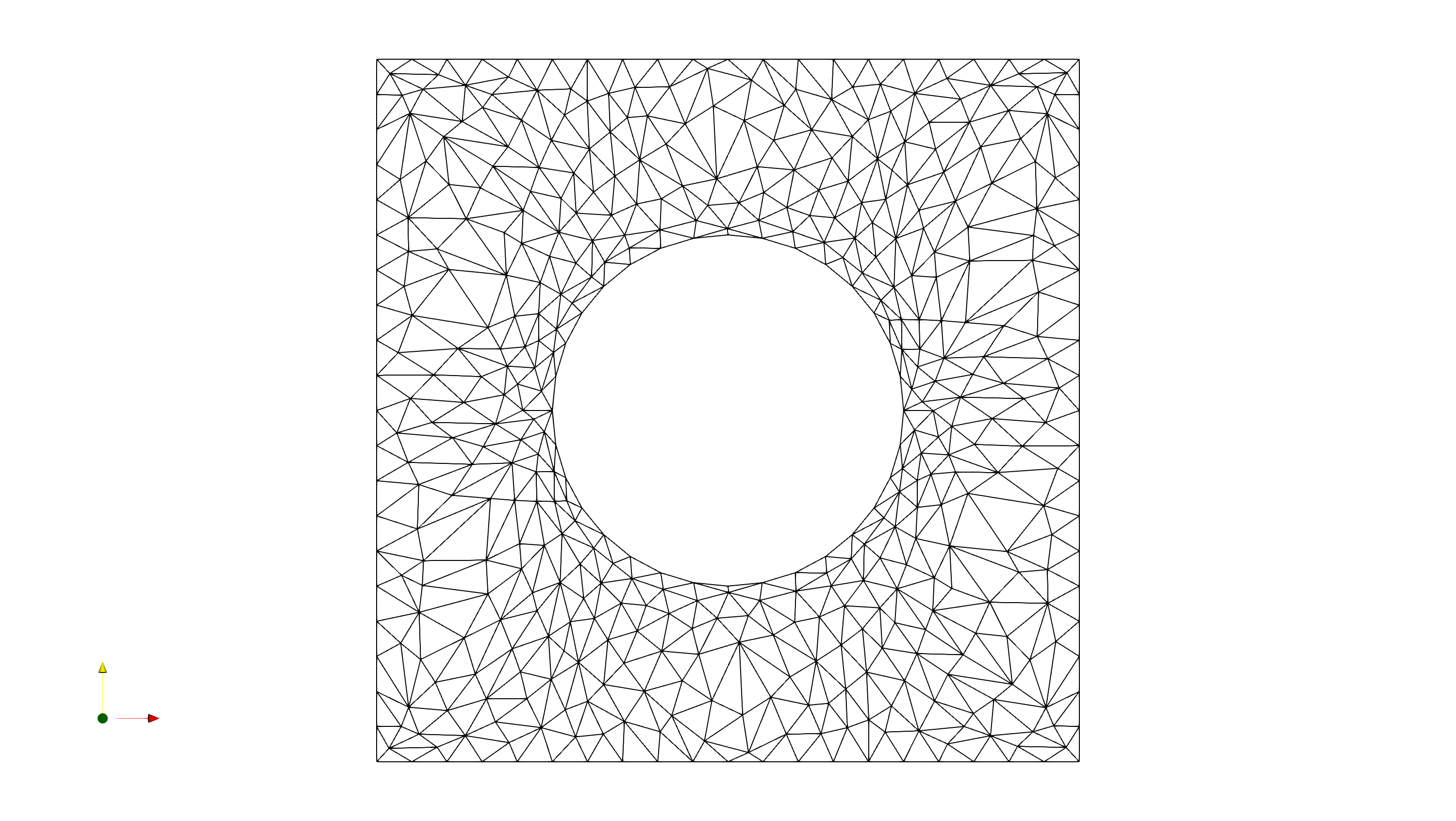} &
			\includegraphics[width=0.3\linewidth, trim={27cm 2cm 27cm 0cm},clip]{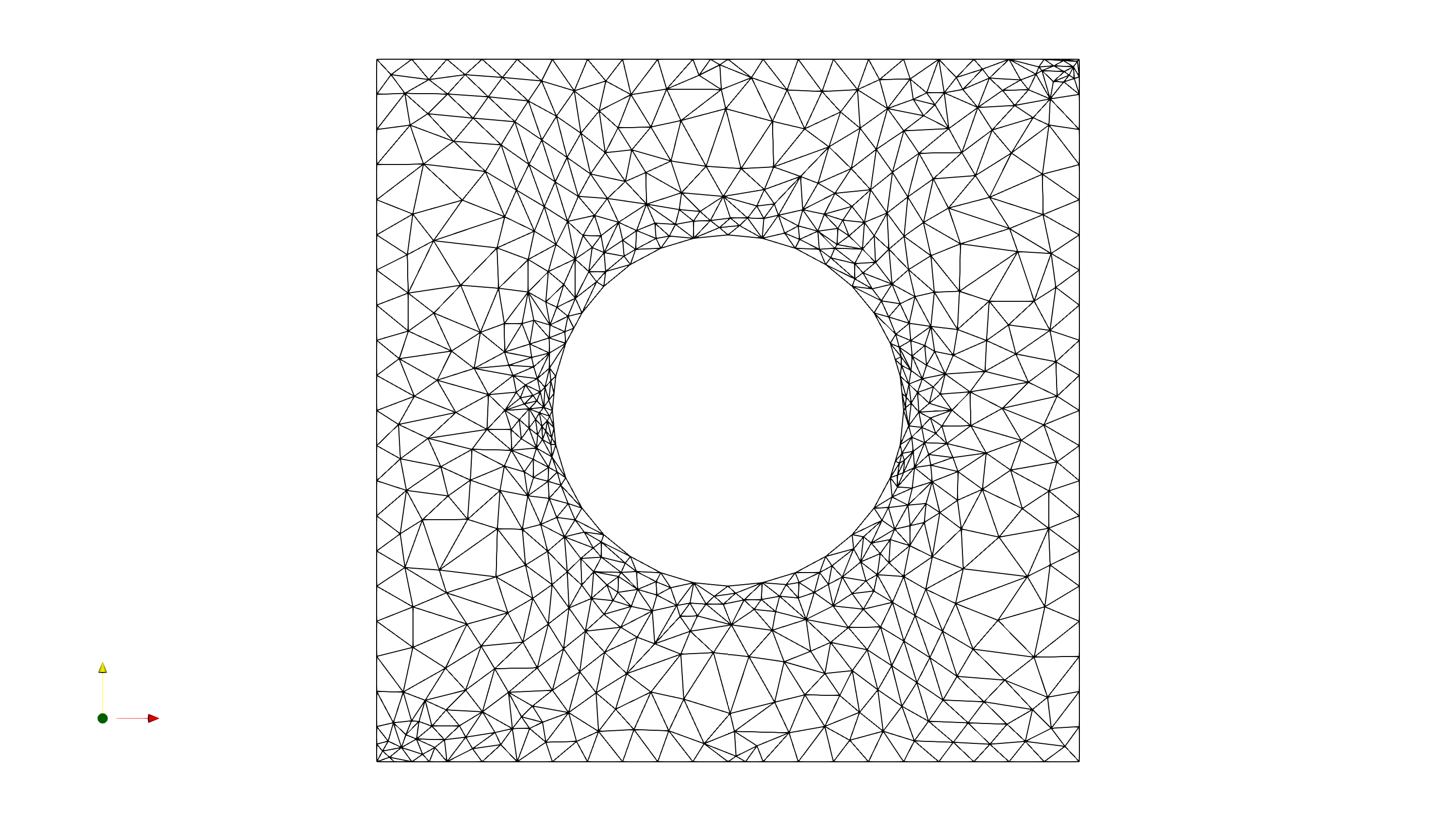} \\
			1 804 nodes & 1 804 nodes & 2 176 nodes \\\hline
			\includegraphics[width=0.3\linewidth, trim={27cm 2cm 27cm 0cm},clip]{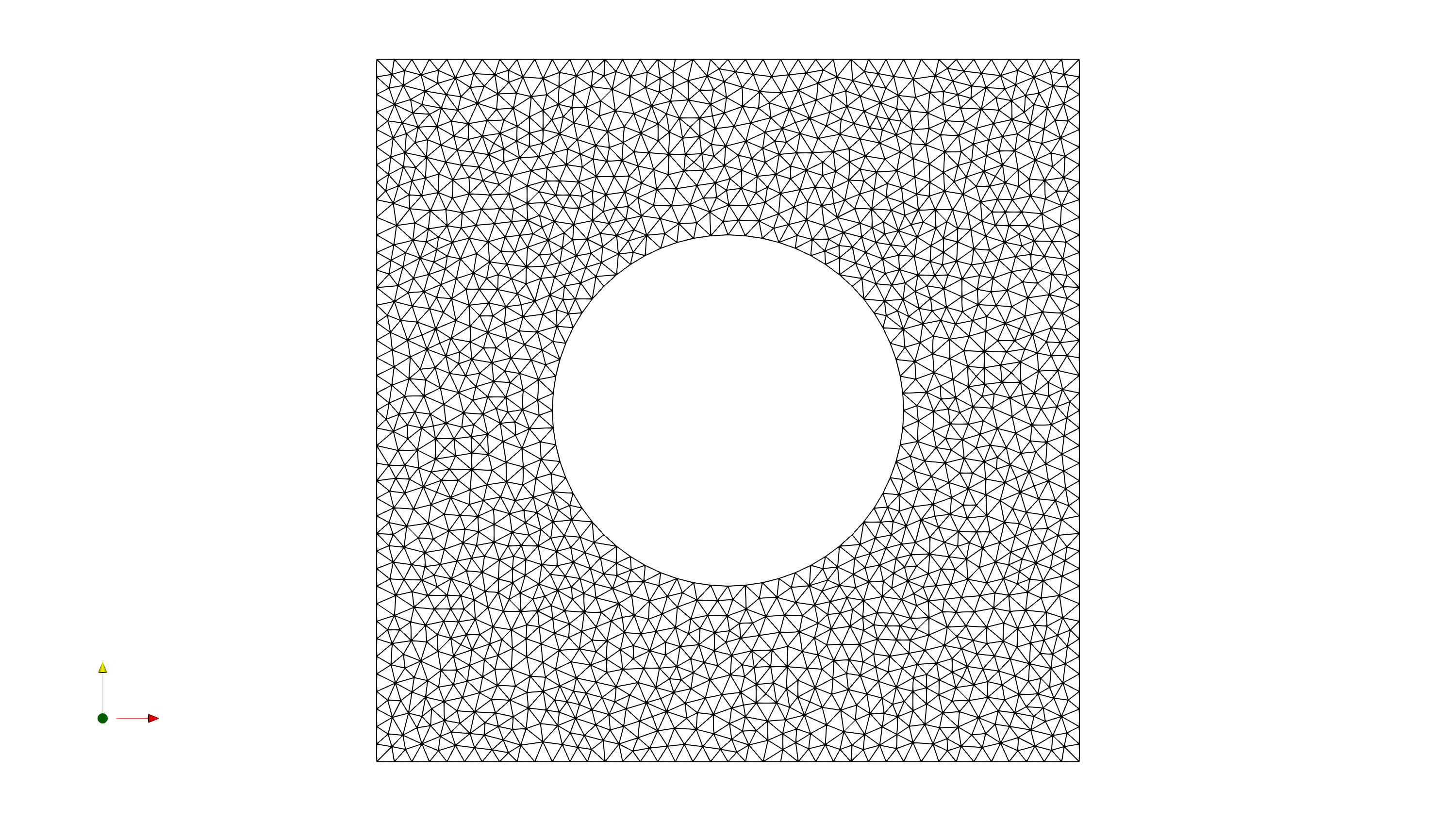} &
			\includegraphics[width=0.3\linewidth, trim={27cm 2cm 27cm 0cm},clip]{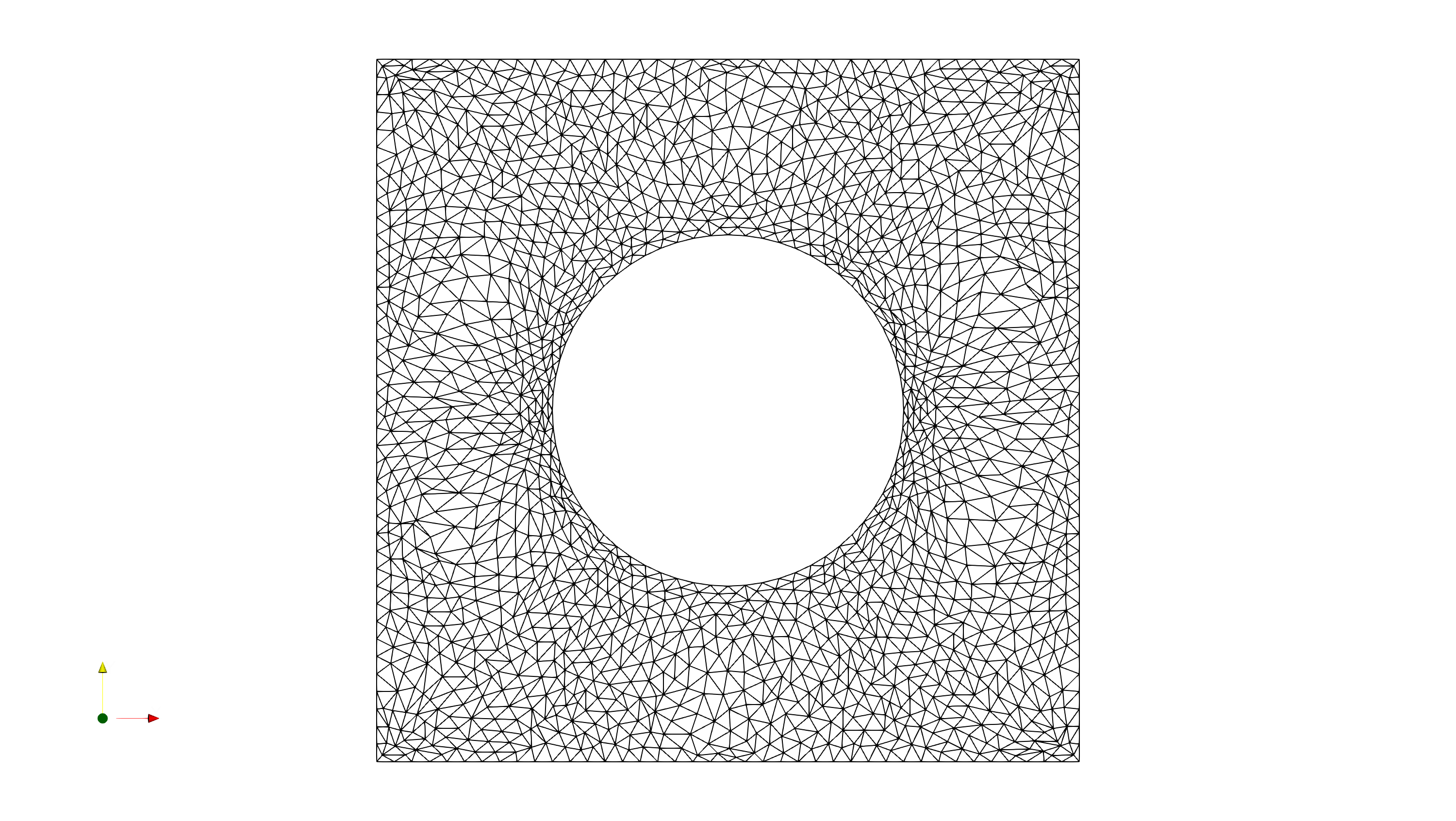} &
			\includegraphics[width=0.3\linewidth, trim={27cm 2cm 27cm 0cm},clip]{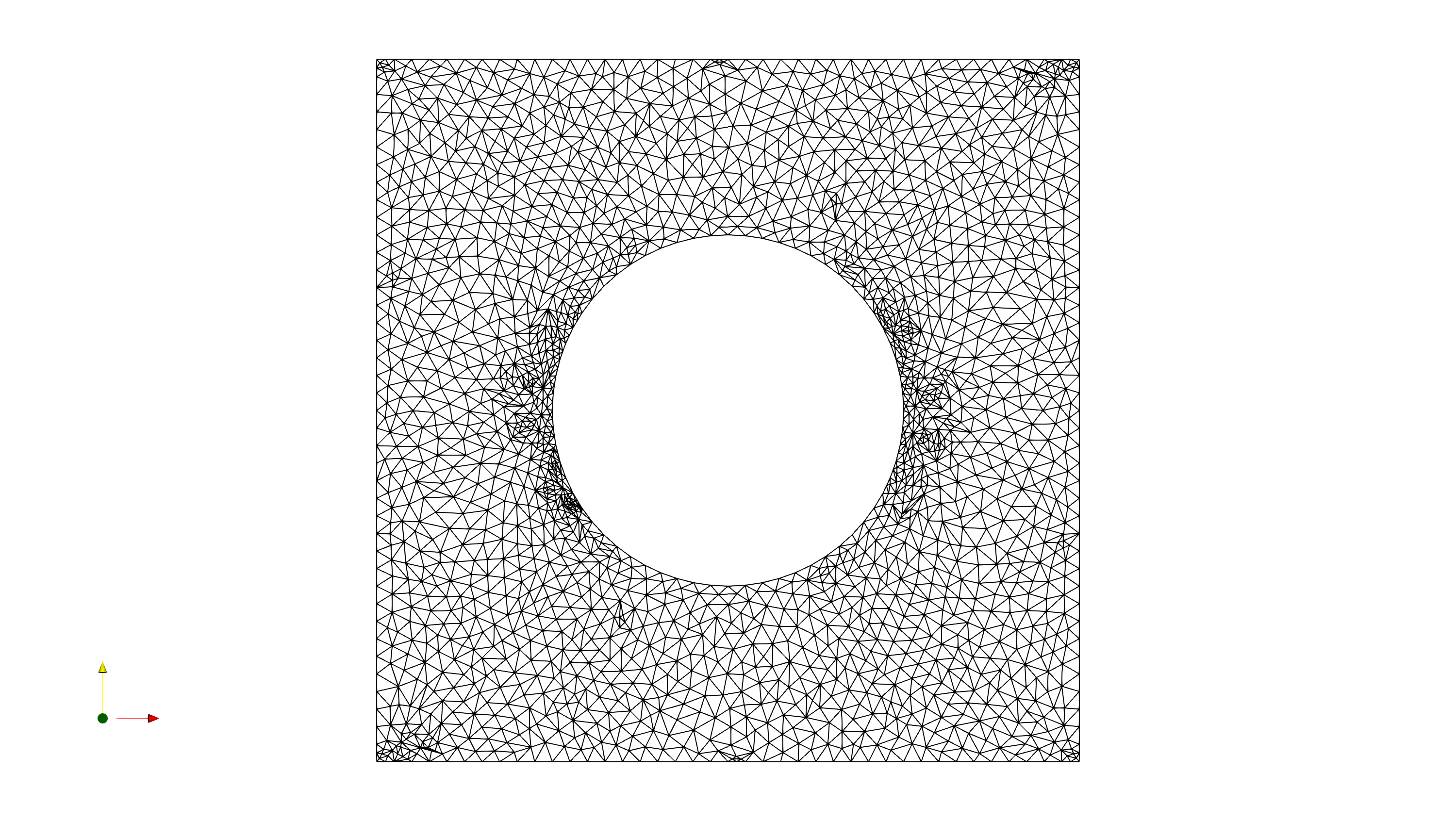} \\
			7 088 nodes & 7 088 nodes & 7 205 nodes \\\hline
			\includegraphics[width=0.3\linewidth, trim={27cm 2cm 27cm 0cm},clip]{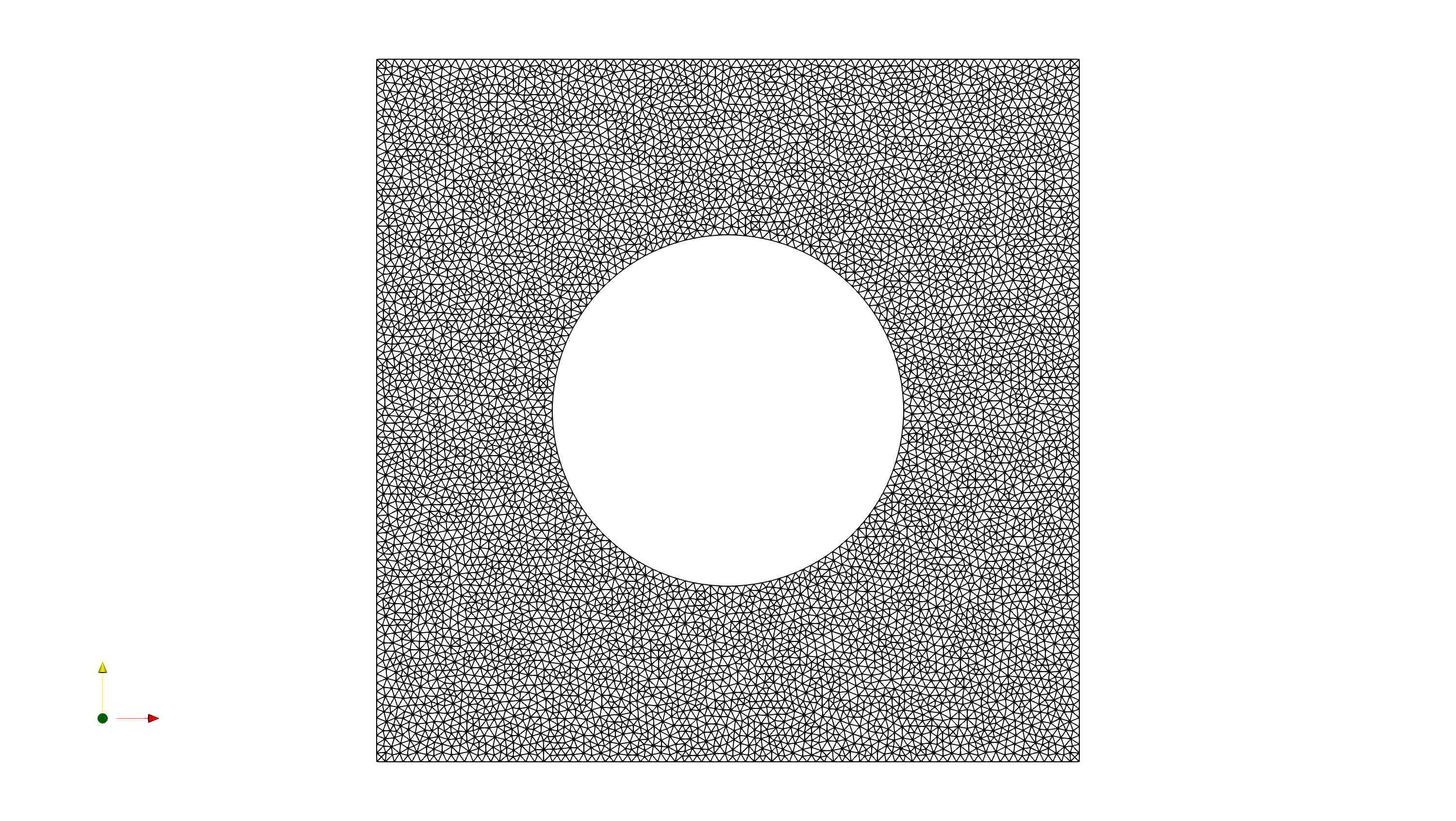} &
			\includegraphics[width=0.3\linewidth, trim={27cm 2cm 27cm 0cm},clip]{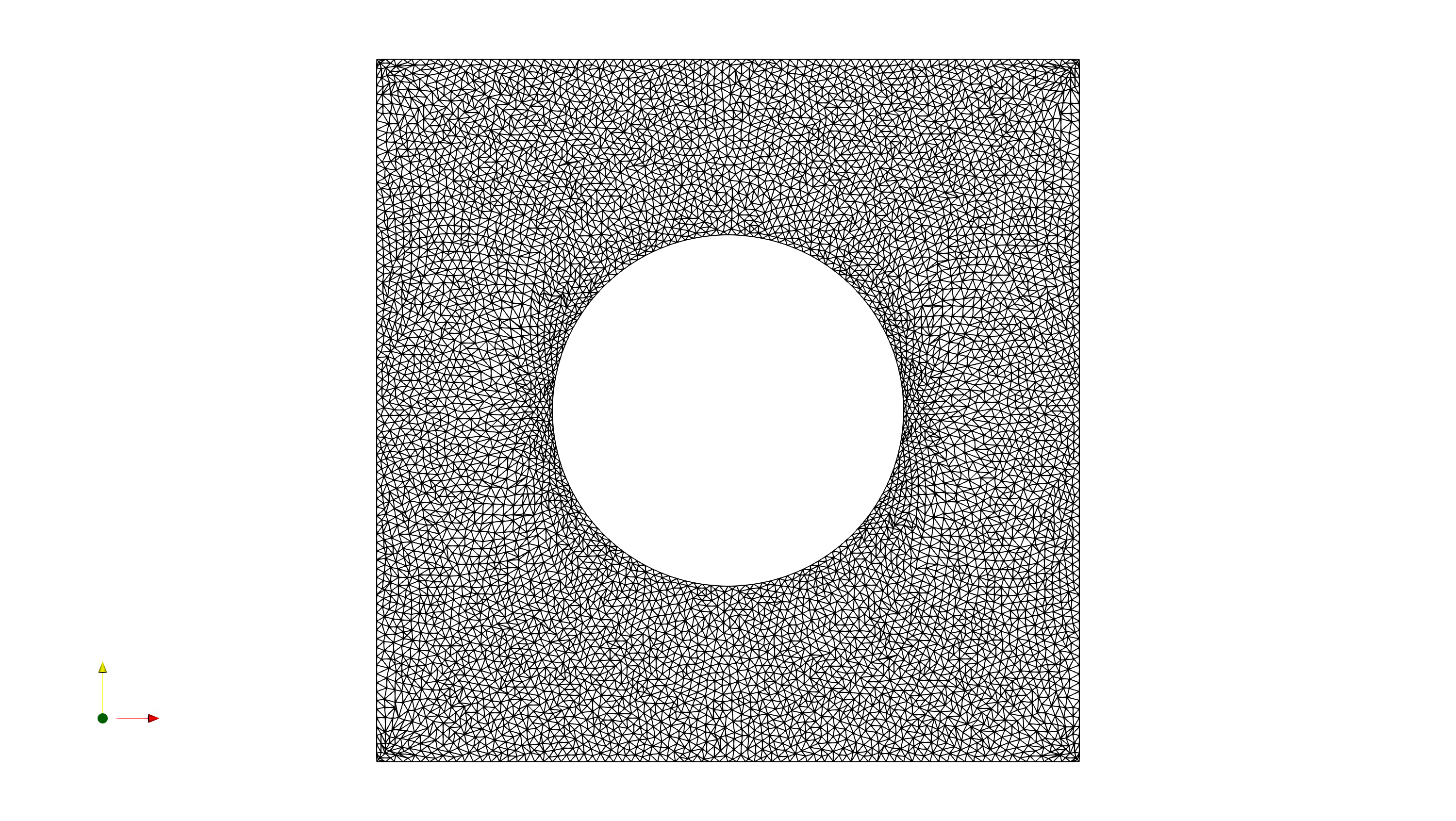} &
			\includegraphics[width=0.3\linewidth, trim={27cm 2cm 27cm 0cm},clip]{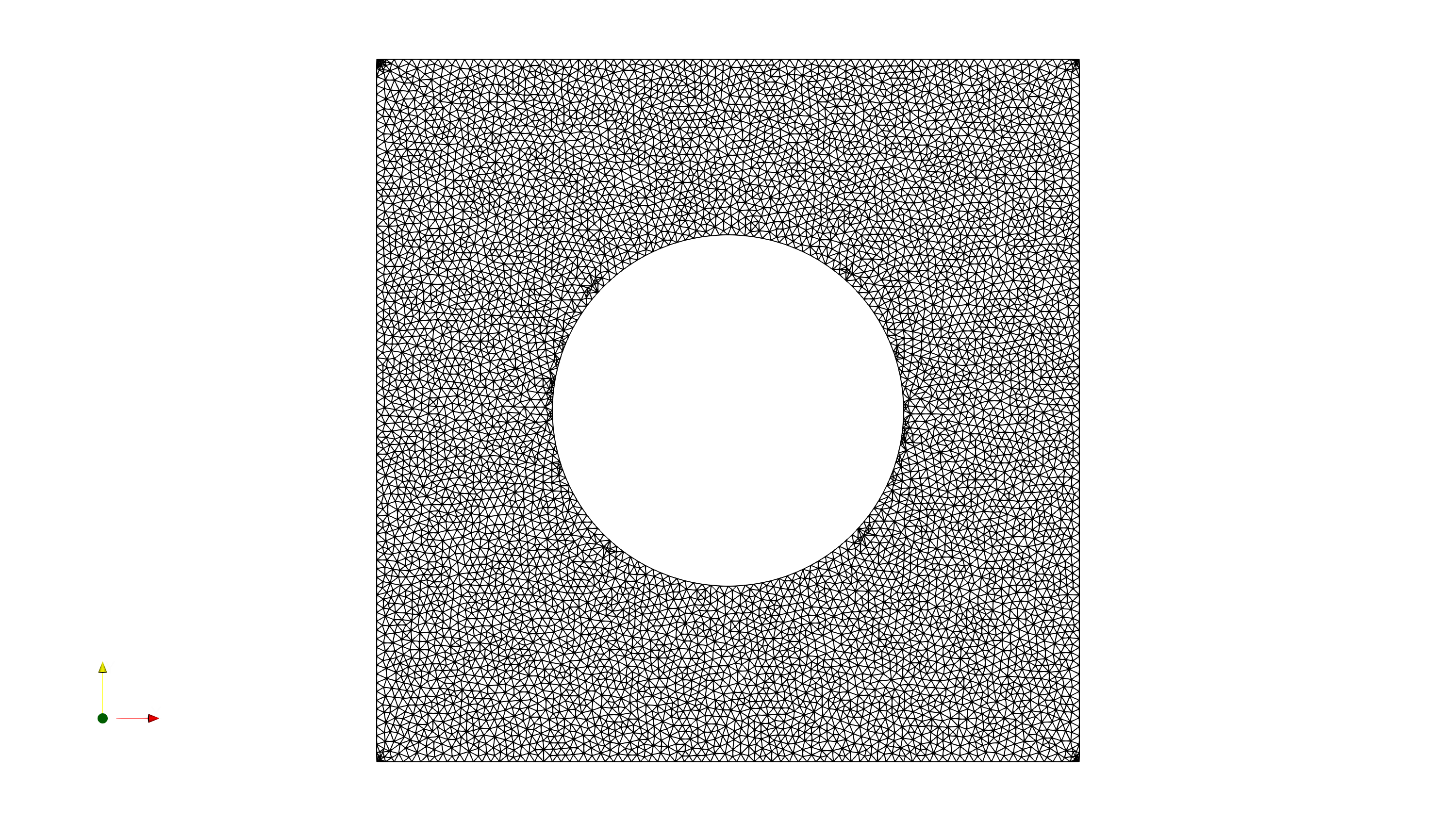} \\
		\end{tabular}
		\captionof{table}{The original discretization mesh of domain $\Omega$ and the final mesh obtained with r- and rh-adaptivity.}
		\label{tab:2D_rh_adapt}
	\end{minipage}
\end{figure}

\KS{
	\subsection{Impact of rh-mesh adaptivity }
	In this section we present the effect of the rh-adaptivity on the 2D problem. The training was performed using the multigrid training strategy and L-BFGS optimizer. The maximal number of element splits was set to 1.
}

The rh-adaptivity led to more accurate results, compared to r-adaptivity, for all except the finest mesh resolution, as shown in Figure \ref{fig:2D_rh_norms}. 
As shown in Table \ref{tab:2D_rh_adapt}, the meshes are refined at the regions where the stress is the highest. For the coarsest mesh resolution, the number of nodes in the rh-adaptive mesh is approximately 2 times higher than in the original mesh. For the finer meshes, the number of added nodes decreases.

\KS{Table \ref{tab:2D_CPU_times} compares the computational time required to obtain a solution on a fixed mesh using a classical FEM solver with the training time for the FENNI models. The comparison shows that, as the number of DOFs increases, the computational time for the FEM solver grows significantly more slowly than the training time of FENNI model. The table also includes training times for adaptive mesh models when employing the multigrid training strategy. Even with the multigrid training, simultaneous optimizing of nodal coordinates and values results in a substantial increase in training time. However, given the accuracy of the model prediction, the training time remains acceptable.}

\KS{We note that while training time could be further reduced by using GPUs or adopting more advanced transfer learning strategies, the real potential of this framework lies in its future application to parametric PDEs. In such cases, a trained model could be used to solve PDEs with different parameters, making the training time more affordable compared to rerunning the FEM solver for each new scenario.}

\begin{table}[h!]
	\centering
	\KS{
		\begin{tabular}{ c c c c | c c c | c c} 
			\multicolumn{4}{c|}{Fixed mesh} & \multicolumn{3}{c|}{r-adaptive mesh} & \multicolumn{2}{c}{rh-adaptive mesh}\\
			DOFs	&	FEniCS	&	Adam	&	L-BFGS	&	DOFs	&	Adam	&	L-BFGS	&	DOFs & L-BFGS	\\
			88	&	0.014	&	2.659	&	0.037	&	88	&	16.63	&	0.222	&	148	&	0.665	\\
			288	&	0.015	&	2.531	&	0.093	&	288	&	20.683	&	4.508	&	584	&	8.930	\\
			968	&	0.018	&	3.170	&	0.228	&	968	&	23.455	&	9.000	&	1368	&	15.226	\\
			3608	&	0.032	&	5.199	&	0.720	&	3608	&	25.381	&	28.11	&	4352	&	54.864	\\
			14176	&	0.083	&	12.155	&	3.067	&	14176	&	27.302	&	76.132	&	14410	&	67.579	\\
		\end{tabular}
		\caption{Comparison of computational times of FEM solver (FEniCS) and training times [s]. Computations were done on CPU (Apple M2, double precision).}
		\label{tab:2D_CPU_times}
	}
\end{table}

\section{Discussion}

The presented numerical experiments show that the choice of loss function and optimizer can have a significant impact on the process of the training as well as the accuracy of results. From the tested settings, the best results on both the fixed and adaptive meshes were obtained when using the combination of potential energy loss and the L-BFGS optimizer. 

The use of weak formulation loss leads to significantly slower convergence. However, the accuracy of the results is comparable to the potential energy loss. The use of the residual loss function was shown to be more problematic\AD{, which is explained by the fact that it contains both convergence and discretization errors}.
Thanks to the model architecture, the Dirichlet boundary conditions do not need to be included in the loss function, as is typical for PINNs. The loss function, therefore, consists only of the PDE residual and compatibility term. Balancing these two terms, however, remains problematic, resulting in different accuracy of predicted displacement and stress. 
We note that only one set of weights was tested in this paper, and thus, a different set of empirically found weights may lead to better results. However, the effort required to find the optimal setting for a given problem must be taken into account. The problem could also be addressed by applying some of the more advanced balancing strategies proposed for PINNs.

The experiments of different mesh resolutions showed, that with the increasing complexity of the problem, more complex training strategies become necessary. In the 2D problem with r- and rh-adaptivity, the convergence for fine meshes was reached only when using the multigrid strategy. Thanks to the interpretability model, other transfer learning strategies could be implemented to speed up the training in the future. \AD{The combined rh-adaptivity and multigrid strategy amounts to optimizing the architecture of the model during the training stage. Such dynamics optimization is a direct consequence of the interpretability offered by the EFENN framework, which enables relying on great transfer learning capabilities, thus overcoming a major limitation that usually comes with standard PINNs.}

In this work, the results are provided for 1D and 2D linear elasticity problems. To solve non-linear problems, the architecture of the model would not be changed as only the loss function would need to be adapted. More generally, to be applied in real-world applications, the incorporation of parameters of the governing equations needs to be addressed. \KS{We showed that the computational cost would be significantly higher compared to using a classical FEM solver if the training needed to be repeated for each set of parameters. This challenge is addressed in Part II  \parencite{daby-seesaramFiniteElementNeural2024}.  }
In \parencite{zhang2022hidenn}, the authors propose to use a fixed tensor decomposition within the HiDeNN framework for the space variables only, allowing the use of separate HiDeNN models for each dimension on meshes generated by the Cartesian product. \AD{In Part II \parencite{daby-seesaramFiniteElementNeural2024}, we propose to extend the concept to the Proper Generalised decomposition (PGD) for surrogate modelling.} The tensor decomposition is greedily built on the space-parameters space, combining separate models for space coordinates and each of the parameters in a versatile framework that builds the reduced-order basis on the fly.

\section{Conclusion}

We presented a machine learning-based method for PDE solutions. The method combines the strengths of two approaches - physics-informed machine learning and FEM.
The framework is based upon a concept of the HiDeNN model \parencite{zhang2021hierarchical}. The model can be seen as a sparse neural network, where the parameters are directly the nodal values and coordinates. Therefore, the training of the model corresponds to simultaneously finding a PDE solution on the optimized mesh.
The concept of HiDeNN is extended in several ways. The original architecture, relying on the construction of global shape functions, was replaced by reference-element-based architecture. Thanks to the implementation of the reference element, the \AD{integral in the }loss function can be easily evaluated using the Gaussian quadrature rule. This was shown to be beneficial, especially when training the model with adaptive mesh. 
Thanks to the new architecture, the model could be straightforwardly extended to 2D unstructured triangular meshes. 
Notably, we \AD{extended} the concept of combined rh-adaptivity \AD{in} 2D \AD{ by designing a r-adaptivity-based refinement criterion.} \AD{Finally, we proposed a rational multigrid training strategy exploiting the interpretable nature of the FENNI}.
The original and new architectures were compared in 1D experiments, showing that the proposed architecture combined with Gaussian quadrature integration provides superior results, especially for adaptive mesh. 
The 2D experiments showed the benefits of the multigrid training strategy in training with adaptive meshes. 
Finally, we discuss how the model could be used for solving parametric PDEs, which would make it significantly more applicable to real-world problems. 

\section*{Acknowledgements}

This work was supported by the French National Research Agency (ANR) under the grant MLQ-CT (ANR-23-CE17-0046) and Bertip EUR (ANR 18EURE0002). 

	\printbibliography[heading=bibintoc]
	
\end{document}